\documentclass{article}
\usepackage{E2VEM}
\usepackage[margin=1in]{geometry}

\bibliography{bibliografia}
\title{SUPG-stabilized stabilization-free VEM: a numerical investigation}
\author{Andrea Borio, Martina Busetto and Francesca Marcon}
\date{}

\begin{document}

\maketitle

\begin{abstract}
  We numerically investigate the possibility of defining stabilization-free
  Virtual Element (VEM) discretizations of advection-diffusion problems in the
  advection-dominated regime. To this end, we consider a SUPG stabilized
  formulation of the scheme. Numerical tests comparing the proposed method with
  standard VEM show that the lack of an additional arbitrary stabilization term, typical of VEM schemes, that adds artificial
  diffusion to the discrete solution, allows to better approximate
  boundary layers, in particular in the case of a low order scheme.
\end{abstract}

\section{Introduction}
\label{sec:intro}
The development of numerical methods for the solution of partial differential
equations exploiting general polygonal or polyhedral meshes has been a topic of
great interest in later years. Among the many families of methods developed in
this context \cite{DIPIETRO2015a,Cicuttin2021,MFD_Brezzi,Riviere2008,PFEM_2004}, 
this paper considers the family of Virtual Element Methods (VEM).

Since the first seminal papers
\cite{Beirao2013a,Beirao2013b,Beirao2014,Beirao2015b}, VEM have been applied in
many contexts where polytopal meshes can be exploited in order to better handle
geometrical complexities of the computational domain, for instance we give a brief, not exhaustive, list of papers \cite{Beirao2015a,Artioli2017,Dassi2020b,Dassi2021b,Benedetto2016b,Benedetto2016c,Benedetto2017,Berrone2023,Borio2011,Berrone2022b}. 
Virtual element schemes are based on the definition of locally computable polynomial projections involved in the discrete bilinear forms. 
These forms consist of the sum of two terms: a singular one, consistent on polynomials, and a stabilizing one that ensures coercivity.
In the literature, the arbitrary nature of the stabilization term
remains an issue to be investigated, since it has been shown that it can cause
problems in many theoretical and numerical contexts. 
For instance, the isotropic nature of the stabilization can become an issue in problems when devising SUPG stabilizations \cite{BBBPSsupg,BBM}.

The aim of this paper is to numerically investigate a possible solution to overcome this issue, conceiving a Virtual Element SUPG formulation that defines coercive bilinear
forms without introducing an arbitrary non-polynomial stabilizing term that in the standard VEM formulation ensures the coercivity of the diffusive part of the problem. 
In the context of the theoretical development of a VEM scheme that do
not require a stabilization term, Stabilization Free Virtual Elements Methods
(SFVEM) have been recently introduced in \cite{BBME2VEM}, in the framework of the lowest order primal discretization of the Poisson equation.
The key idea of the proposed method is to define self-stabilized bilinear forms, exploiting only higher order polynomial projections.
The theoretical study of this method is an ongoing investigation, however some
tests on highly anisotropic problems show that this new approach is able to overcome some issues of standard VEM related to the isotropic nature of the
stabilization \cite{BBME2VEMLetter,BBMTLetter}.

The outline of the paper is as follows. In Section \ref{sec:modelproblem} we present the advection-diffusion model problem. In Section \ref{sec:discretization} we define the numerical scheme. In Section \ref{sec:wellpos} we discuss the weel-posedness of the discrete problem. 
Section \ref{sec:apriori} is devoted to a priori error estimates. Finally, in
Section \ref{sec:numres} we present some numerical results that assess the advantages of the presented method.

In the following, $\scal[\omega]{\cdot}{\cdot}$ denotes the
$\lebl{\omega}$-scalar product, $\norm[\omega]{\cdot}$ denotes the corresponding
norm, and $\norm[m,\omega]{\cdot}$ and $\seminorm[m,\omega]{\cdot}$ denote the
$\sobh{m}{\omega}$ norm and semi-norm.

\section{Model Problem}
\label{sec:modelproblem}

Let $\Omega \subset \mathbb{R}^2$ be a bounded open set. We consider the
following advection-diffusion model problem: find $u$ such that
\begin{equation}
  \label{eq:theproblem}
  \begin{cases}
    -\K\Delta u + \beta\cdot \nabla u = f & \text{in $\Omega$\,,}
    \\
    u = 0 & \text{on $\partial \Omega$\,,}
  \end{cases}
\end{equation}
where $\K > 0$ is a positive real number and we make the standard hypotheses that $\beta \in \left[\lebl[\infty]{\Omega}\right]^2$,
$\div\beta = 0$, and $f\in\lebl{\Omega}$. 
Moreover, we define the bilinear form
$\B{}{}\colon \sobh[0]{1}{\Omega} \times \sobh[0]{1}{\Omega} \to \mathbb{R}$ and
the operator $\F{}\colon \lebl{\Omega} \to \mathbb{R}$ such that
\begin{align*}
  \B{w}{v} &= \scal[\Omega]{\K \nabla w}{\nabla v} + \scal[\Omega]{\beta\cdot\nabla w}{v} \quad
             \forall w,v \in \sobh[0]{1}{\Omega} \,,
  \\
  \F{v} &= \scal[\Omega]{f}{v} \quad \forall v \in \lebl{\Omega} \,.
\end{align*}
Then, the variational formulation of \eqref{eq:theproblem} reads as follows:
find $u\in \sobh[0]{1}{\Omega}$ such that
\begin{equation}
  \label{eq:modelvarform}
  \B{u}{v} = \F{v} \quad \forall v \in \sobh[0]{1}{\Omega} \,.
\end{equation}
It is a standard result that the above problem is well-posed under the above
regularity assumptions on the data. For the sake of better readibility, we limit ourselves to homogeneous
Dirichlet boundary conditions, but more
general boundary conditions can be considered and will be considered in the
numerical tests.  Finally, for any open subset $\omega \subset \Omega$, we
define:
\begin{equation*}
  \beta_\omega = \sup_{v \in \lebldouble{\omega}} \frac{\norm[\omega]{\beta \cdot v}}{\norm[\omega]{v}} \,.
\end{equation*}

\section{Problem discretization}
\label{sec:discretization}

This section is devoted to the discretization of \eqref{eq:modelvarform} using
an enlarged enhancement Virtual Element space. The discretization
defined here was introduced in its lowest order version in
\cite{BBME2VEM}. Here, we generalize that scheme to a generic order $k \geq 1$.

\subsection{Discrete space}
\label{sec:E2VEMspace}

Let $\Mh$ denote a mesh of $\Omega$ made up of polygons. Let $h_E$ denote the
diameter of $E\in\Mh$, and let $h = \max_{E\in\Mh} h_E$ be the mesh
parameter. We make the following mesh assumptions, that are standard in the
context of VEM \cite{Beirao2013a,Beirao2015b}: there exists a constant
$\kappa > 0$ independent of $h$ such that, for any polygon $E\in\Mh$, if
$\Eh[E]$ denotes the set of edges of $E$,
\begin{enumerate}
\item $E$ is star-shaped with respect to a ball of radius
  $\rho \geq \kappa h_E$;
\item $\forall e \in\Eh[E]$, $\abs{e} \geq \kappa h_E$, where $\abs{e}$ denotes
  the length of $e$.
\end{enumerate}

For each $E\in\Mh$, let $\Poly{n}{E}$ be the space of polynomials of degree at
most $n$, for any $n\in\mathbb{N}$. As a basis of $\Poly{n}{E}$, we
choose the following set of scaled monomials:
\begin{equation*}
  \Monom{n}{E} = \left\{m_{\bs\alpha}\colon \mathbb{R}^2\to \mathbb{R} \text{ such that }
    m_{\bs\alpha}(x,y) = \frac{(x-x_E)^{\alpha_1}(y-y_E)^{\alpha_2}}{h_E}
 \vphantom{\frac{(x-x_E)^{\alpha_1}(y-y_E)^{\alpha_2}}{h_E}},\,\,
 0 \leq \abs{\bs\alpha} = \alpha_1+\alpha_2 \leq n \right\} \,,
\end{equation*}
where $(x_E,y_E)$ is the center with respect to which $E$ is star-shaped.
Moreover, let us define a subspace of $\Monom{n}{E}$ given by
\begin{equation*}
  \Monom{n,t}{E} = \left\{m_{\bs\alpha} \in \Monom{n}{E}\colon
    \abs{\bs\alpha} \geq t \right\} \quad \forall t<n\,.
\end{equation*}
Let $\proj[\nabla,E]{n}{}\colon \sobh{1}{E} \to \Poly{n}{E}$ be the projection
operator such that, $\forall v \in \sobh{1}{E}$,
\begin{equation}
  \label{eq:defPiNabla}
  \begin{cases}
    \scal[E]{\nabla \proj[\nabla,E]{n}{}v}{\nabla p} = \scal[E]{\nabla v}{\nabla p}
    & \forall p\in\Poly{n}{E}\,,
    \\
    \int_{\partial E}\proj[\nabla,E]{n}{}v = \int_{\partial E} v & \text{if $k = 1$} \,,
    \\
    \int_{E}\proj[\nabla,E]{n}{}v  = \int_{E} v & \text{if $k > 1$} \,.
  \end{cases}
\end{equation}
 Then, given
$\ell_E\in\mathbb{N}$, $\ell_E\geq 0 $, we introduce the enlarged enhancement VEM discrete
space:
\begin{multline}
  \label{eq:localE2VEMspace}
  \Vh[E]{k,\ell_E} =
  \left\{
    v \in \sobh{1}{E}\colon \Delta v \in\Poly{k+\ell_E}{E}\,,
  \trace[\partial E]{v}\in\cont{\partial E}\,,
  \trace[e]{v}\in\Poly{k}{e} \;\forall e \in\Eh[E]\,, \right.
  \\
  \left. \scal[E]{v}{p} = \scal[E]{\proj[\nabla,E]{k}{}v }{p}
    \; \forall p \in\Poly{k+\ell_E,k-2}{E}\right\} \,,
\end{multline}
where $\gamma$ denotes the trace operator and
\begin{equation*}
  \Poly{k+\ell_E,k-2}{E} = \Span\Monom{k+\ell_E,k-2}{E} \,.
\end{equation*}

\begin{remark}
  The space $\Vh[E]{k,0}$ is the standard VEM space used in \cite{Beirao2015b,BBBPSsupg}.
\end{remark}

Notice that, $\forall \ell_E$, a possible set of degrees of freedom of
$\Vh[E]{k,\ell_E}$ corresponds to the one introduced for standard VEM spaces
(see e.g. \cite{Beirao2015b,BBBPSsupg}), that is, for any $v_h\in\Vh[E]{k,\ell_E}$,
\begin{enumerate}
\item the values of $v$ at the vertices of $E$;
\item for each $e\in\Eh[E]$, the values of $v$ at $k-1$ points internal to $e$;
\item the scaled moments $\frac{1}{\abs{E}}\scal[E]{v}{m_{\bs\alpha}}$
  $\forall m_{\bs\alpha}\in\Monom{k-2}{E}$.
\end{enumerate}

\subsection{Discrete problem}
\label{sec:discrprob}

This section is devoted to define the SUPG-stabilized discrete version of
\eqref{eq:modelvarform}. For this purpose, we adapt the scheme proposed and
analysed in \cite{BBBPSsupg}. Notice that the discretization presented here
differs from the cited one only in the bilinear form discretizing the diffusive
part of the problem, while the other components of the discrete bilinear form
are untouched, and so is the discrete right-hand side.

Let $E\in \Mh$ be given. If $k>1$, let $\tilde{C}_{k}$ be the largest constant, independent of $h_E$, 
such that
\begin{equation}
  \label{eq:inverseinequality}
  \tilde{C}_{k}h^2_E \norm[E]{\Delta p}^2 \leq
  \norm[E]{\nabla p}^2 \quad \forall p \in \Poly{k}{E}\,.
\end{equation}
Such an inequality is a standard inverse inequalities for polynomials (see
e.g. \cite{BrennerScott}). In particular, $\tilde{C}_{k}$ depends on the
constant $\kappa$ in the mesh assumptions, and on the degree $k$. Notice that it
does not depend on $\ell_E$. The \emph{P\'eclet} number associated to $E$ is
defined as
\begin{equation*}
  \Pe[E] = m_k \frac{\beta_E h_E}{\K} \,, \quad
  m_k =
  \begin{cases}
    \frac13 & \text{if $k=1$,} \\ 2\tilde{C}_k & \text{if $k>1$.}
  \end{cases}
\end{equation*}
Following the standard Streamline Upwind Petrov-Galerkin approach
\cite{franca1992stabilized}, for any $E\in\Mh$ let the parameter $\tau_E$ be
given, such that
\begin{equation}
  \label{eq:deftau}
  \tau_E = \frac{h_E}{2\beta_E} \min\left\{1,\Pe[E]\right\} \,,
\end{equation}
and we define the symmetric positive definite tensor
$\K[\beta,E] = \K + \tau_E\beta\beta^\intercal$. Notice that, by the hypotheses
on $\beta$, this tensor satisfies, $\forall E \in \Mh$,
$\forall \bs{v} \in \lebldouble{E}$,
\begin{align*}
   \sqrt{\K}\norm[E]{\bs{v}} \leq \norm[E]{\sqrt{\K[\beta,E]} \bs{v}}
  \leq \sqrt{\left(\K + \tau_E \beta_E^2\right)} \norm[E]{\bs{v}}\,.
\end{align*}
To define the SUPG formulation of the problem, we introduce the space
\begin{equation*}
  \sobhDeltaloc{\Mh} = \left\{v \in \sobh[0]{1}{\Omega} \colon
    \Delta v \in \lebl{E} \quad \forall E \in \Mh\right\} \,,
\end{equation*}
and, $\forall E \in \Mh$, the bilinear form
$\Bsupg[E]{}{}\colon \sobhDeltaloc{\Mh}\times \sobh[0]{1}{\Omega} \to
\mathbb{R}$ such that, $\forall w \in \sobhDeltaloc{\Mh}$,
$\forall v \in \sobh[0]{1}{\Omega}$,
\begin{equation*}
  \Bsupg[E]{w}{v} = \a[E]{w}{v} + \b[E]{w}{v} + \d[E]{w}{v} \,,
\end{equation*}
where
\begin{align*}
  \a[E]{w}{v} &= \scal[E]{\K\nabla w}{\nabla v} +
                \tau_E \scal[E]{\beta \cdot \nabla w}{\beta \cdot \nabla v} \,,
  \\
  \b[E]{w}{v} &= \scal[E]{\beta \cdot \nabla w}{v} \,,
  \\
  \d[E]{w}{v} &= -\tau_E\scal[E]{\K\Delta w}{\beta\cdot\nabla v} \,.
\end{align*}
Moreover, we introduce the operator
$\Fsupg[E]{}\colon \sobh{1}{E} \to \mathbb{R}$ such that
\begin{equation*}
  \Fsupg[E]{v} = \scal[E]{f}{v + \tau_E\beta\cdot\nabla v} \quad
  \forall v \in \sobh{1}{E} \,.
\end{equation*}
The above operators are not computable from the degrees of freedom for functions in $\Vh[E]{k,\ell_E}$. For
this reason, we define discrete bilinear forms involving polynomial projections
of VEM discrete functions. In the spirit of \cite{BBME2VEM}, for a given $n\in\mathbb{N}$ let
$\proj[0,E]{n}{}\colon \lebldouble{E} \to \Polydouble{n}{E}$ denote the
$\lebldouble{E}$ projection onto $\Polydouble{n}{E}$ and let
\begin{align*}
  \ah[E]{w}{v} &= \scal[E]{\K\proj[0,E]{k+\ell_E-1}{}\nabla w}{\proj[0,E]{k+\ell_E-1}{}\nabla v} +
                 \tau_E \scal[E]{\beta \cdot \proj[0,E]{k+\ell_E-1}{}\nabla w}{\beta \cdot \proj[0,E]{k+\ell_E-1}{}\nabla v} \,,
  \\
  \bh[E]{w}{v} &= \scal[E]{\beta \cdot \proj[0,E]{k-1}{}\nabla w}{\proj[0,E]{k-1}{} v} \,,
  \\
  \dh[E]{w}{v} &= -\tau_E\scal[E]{\K\div\left(\proj[0,E]{k-1}{}\nabla w\right)}{\beta\cdot\proj[0,E]{k+\ell_E-1}{}\nabla v} \,.
\end{align*}
\begin{remark}
  Notice that, contrarily to standard VEM SUPG formulations
  (cf. \cite{BBBPSsupg,BBM}), here we do not require an additional arbitrary stabilization
  term in the bilinear form $\ah[E]{}{}$. The choice of the parameter $\ell_E$
  is done in order to guarantee coercivity of $\ah[E]{}{}$, as detailed in the
  next section.
\end{remark}
Finally, to state our discrete problem, we define the bilinear form
$\Bsupgh[E]{}{}\colon \sobhDeltaloc{\Mh}\times \sobh[0]{1}{\Omega} \to
\mathbb{R}$ such that, $\forall w \in \sobhDeltaloc{\Mh}$,
$\forall v \in \sobh[0]{1}{\Omega}$,
\begin{equation*}
  \Bsupgh[E]{w}{v} = \ah[E]{w}{v} + \bh[E]{w}{v} + \dh[E]{w}{v} \,,
\end{equation*}
and the operator $\Fsupgh[E]{}\colon \sobh{1}{E} \to \mathbb{R}$ such that
\begin{equation*}
  \Fsupgh[E]{v} = \scal[E]{f}{\proj[0,E]{k-1}{}v +
    \tau_E\beta\cdot\proj[0,E]{k+\ell_E-1}{}\nabla v} \quad
  \forall v \in \sobh{1}{E} \,.
\end{equation*}
Then, let
\begin{equation*}
  \Vh{k,\ell} = \left\{v\in\sobh[0]{1}{\Omega}\colon v \in
    \Vh[E]{k,\ell_E} \; \forall E \in \Mh\right\}\,,
\end{equation*}
our discretization of \eqref{eq:modelvarform} requires to find
$u_h \in \Vh{k,\ell}$ such that
\begin{equation}
  \label{eq:discrvarform}
  \sum_{E\in\Mh}\Bsupgh[E]{u_h}{v_h} = \sum_{E\in\Mh} \Fsupgh[E]{v_h} \quad
  \forall v_h \in \Vh{k,\ell} \,.
\end{equation}


\section{Well-posedness}
\label{sec:wellpos}

The aim of this section is to discuss the key points needed for the well-posedness of \eqref{eq:discrvarform}. 
It is currently an open problem to prove theoretically a robust criterium to choose $\ell_E$ for any kind of polygon, in such a way that $\ah[E]{}{}$ is coercive. 
Thus, we assume that there exists at least a good choice of $\ell_E$ for each type of polygon and in Section \ref{sec:numres} we perform a numerical investigation of it.

\begin{assumption}
  \label{eq:assumcoercivity}
  We assume that, $\forall k \geq 1$ $\exists \ell_E$ such that
  $\exists \alpha > 0$ independent of $h_E$ satisfying
  \begin{equation}
    \label{eq:localcoercivityallk}
    \norm[E]{\proj[0,E]{k+\ell_E-1}{}\nabla v_h}^2 \geq \alpha \norm[E]{\nabla v_h}^2 \quad \forall v_h \in \Vh[E]{k,\ell_E}\,.
  \end{equation}
\end{assumption}

From now on, we
set $\ell_E$ as the smallest integer satisfying Assumption \ref{eq:assumcoercivity}.
With this choice, we define the following VEM-SUPG norm:
\begin{equation*}
  \ennorm{v_h}^2 = \sum_{E\in\Mh} \limits \ahE{v_h}{v_h} \quad \forall v_h\in\Vh{k,\ell} \,.
\end{equation*}
Then, we can prove the following well-posedness result.
\begin{theorem}
  Under Assumption \ref{eq:assumcoercivity}, we have, $\forall E\in \Mh$ and for
  $h$ sufficiently small,
  \begin{equation*}
    \exists C >0 \colon \sum_{E\in\Mh}\limits\Bsupgh[E]{v_h}{v_h} \geq C\ennorm{v_h}^2
    \quad \forall v_h \in \Vh[E]{k,\ell_E}\,.
  \end{equation*}
\end{theorem}
\begin{proof}
  Let $v_h\in\Vh[E]{k,\ell_E}$ be given. First, exploiting the definition of
  $\tau_E$ in \eqref{eq:deftau}, the inverse inequality
  \eqref{eq:inverseinequality} and Young
  inequalities, we get,
  \begin{multline*}
    \abs{\dh[E]{v_h}{v_h}} = \tau_E \abs{\scal[E]{\K\div\left(\proj[0,E]{k-1}{}\nabla v_h\right)}{\beta\cdot\proj[0,E]{k+\ell_E-1}{}\nabla v_h}} 
    \\
    \leq \frac{m_kh^2_E}{4\K} \norm[E]{\K\div\left( \proj[0,E]{k-1}{}\nabla v_h\right)}^2 + \frac{\tau_E}{2}\norm[E]{\beta\cdot\proj[0,E]{k+\ell_E-1}{}\nabla v_h}^2 
    \leq \frac{1}{2\K} \norm[E]{\K \proj[0,E]{k-1}{}\nabla v_h}^2 + \frac{\tau_E}{2}\norm[E]{\beta\cdot\proj[0,E]{k+\ell_E-1}{}\nabla v_h}^2
    \\
    \leq \frac{\K}{2}\norm[E]{\proj[0,E]{k-1}{}\nabla v_h}^2 + \frac{\tau_E}{2}\norm[E]{\beta\cdot\proj[0,E]{k+\ell_E-1}{}\nabla v_h}^2 
    \\
    \leq \frac{\K}{2}
    \norm[E]{ \proj[0,E]{k+\ell_E-1}{}\nabla v_h}^2
    + \frac{\tau_E}{2}\norm[E]{\beta\cdot\proj[0,E]{k+\ell_E-1}{}\nabla v_h}^2
    = \frac{1}{2} \ah[E]{v_h}{v_h}\,.
  \end{multline*}
  Thus, it follows
  \begin{equation}\label{eq:wellpos:ah+dh}
    \ah[E]{v_h}{v_h} + \dh[E]{v_h}{v_h} \geq \frac{1}{2}\ah[E]{v_h}{v_h} = \frac{1}{2}\ennorm[E]{v_h}^2\,.
  \end{equation}
  Moreover, since $\b{v_h}{v_h} = 0$, we get
  \begin{multline*}
    \abs{ \sum_{E\in\Mh}\bh[E]{v_h}{v_h}} = \abs{\sum_{E\in\Mh} \scal[E]{\beta\cdot\proj[0,E]{k-1}{}\nabla v_h}{\proj[0,E]{k-1}{} v_h}}
    = \abs{\sum_{E\in\Mh} \scal[E]{\proj[0,E]{k-1}{\beta\cdot\proj[0,E]{k-1}{}\nabla v_h}}{v_h}}
    \\
    = \abs{\sum_{E\in\Mh} \scal[E]{\proj[0,E]{k-1}{\beta\cdot\proj[0,E]{k-1}{}\nabla v_h}-\beta\cdot\nabla v_h}{v_h}}
    =  \abs{\sum_{E\in\Mh} \scal[E]{\beta\cdot\proj[0,E]{k-1}{}\nabla v_h-\beta\cdot\nabla v_h}{\proj[0,E]{k-1}{}v_h-v_h}}
    \\
    \leq \sum_{E\in\Mh}\norm[E]{\beta\cdot\left(\proj[0,E]{k-1}{}\nabla v_h - \nabla v_h\right)} \norm[E]{\proj[0,E]{k-1}{}v_h - v_h}
    \leq C_{\K\beta} \sum_{E\in\Mh}h_E \ennorm[E]{v_h}^2 \,,
  \end{multline*}
  where $C_{\K\beta}$ depends in particular on the problem data and the
  equivalence constant in \eqref{eq:localcoercivityallk}. Collecting the latest
  estimate and \eqref{eq:wellpos:ah+dh}, we get the thesis:
  \begin{equation*}
    \sum_{E\in\Mh}\limits\Bsupgh[E]{v_h}{v_h} \geq
    \left( \frac{1}{2} - C_{\K\beta}h \right) \ennorm{v_h}^2 \,.
  \end{equation*}
\end{proof}

Notice that the above result is analogous to the one obtained in
\cite{BBBPSsupg} in the case of standard VEM. Notice that other choices are
possible to discretize the transport term (see \cite{Beirao2021}), and they
would lead to a similar well-posedness result.


\section{Error analysis}
\label{sec:apriori}
In this section we address the a priori error analysis of the presented scheme,
under Assumption \ref{eq:assumcoercivity}. The analysis follows the techniques
already used to prove a priori estimates for standard VEM schemes
\cite{Beirao2015b,BBBPSsupg,Beirao2021}. We provide some
details about the interpolation properties of the considered space, since it is
a slight modification of the classical VEM spaces, due to the enlarged
enhancement property.

In the following we prove an interpolation estimate onto the space
$\Vh[E]{k,\ell_E}$ defined by \eqref{eq:localE2VEMspace}. The proof follows the
one in \cite{Cangiani2017} for the interpolation on standard VEM spaces. First
of all we prove an auxiliary result, introducing an inverse inequality in
$\Vh[E]{k,\ell_E}$.

\begin{lemma}
  \label{lem:inverse-inequality}
  Let $E\in\Mh$ and let $w\in\sobh{1}{E}$ such that
  $\Delta w\in\Poly{k+\ell}{E}$, for some chosen $\ell \geq 0$. Then, there
  exists a constant $C_I$ such that
  \begin{equation}
    \label{eq:inverse-estim}
    \norm[E]{\Delta w} \leq C_I h_E^{-1}\norm[E]{\nabla w} \,.
  \end{equation}
  The constant $C_I$ depends on $k$, $\ell_E$ and on the mesh regularity
  parameter $\kappa$.
\end{lemma}
\begin{proof}
  First, let $p\in\Poly{k+\ell}{E}$. Then, let $\psi_E \in\sobh[0]{1}{E}$ be a
  bubble function defined on a regular sub-triangulation of $E$ (see
  \cite{Cangiani2017}), such that $1\geq\psi_E\geq 0$. Then, since
  $\psi_E p\in \sobh[0]{1}{E}$,
  \begin{equation*}
    \norm[\sobh{-1}{E}]{p} =
    \sup_{v\in\sobh[0]{1}{E}}\frac{\scal[E]{p}{v}}{\norm[E]{\nabla v}}
    \geq \frac{\scal[E]{\psi_E}{p^2}}{\norm[E]{\nabla (\psi_Ep)}} \geq
    C_B \frac{\norm[E]{p}^2}{h_E^{-1}\norm[E]{p}} = C_B h_E \norm[E]{p} \,,
  \end{equation*}
  where we use the fact that, since $\psi_E$ is non-negative and bounded,
  $\sqrt{\scal[E]{\psi_E}{p^2}}$ is a norm on $\Poly{k+\ell}{E}$, and thus
  equivalent to $\norm[E]{p}$. Similarly, since $\psi_E\in \sobh[0]{1}{E}$,
  $\norm[E]{\nabla(\psi_Ep)}$ is a norm on $\Poly{k+\ell}{E}$, and standard
  scaling arguments provide the weight $h_E^{-1}$. It follows that $C_B$ depends
  on $k$, $\ell$ and on the shape regularity of $E$.

  Taking $p=\Delta w$ in the above result and applying Green's theorem and a
  Cauchy-Schwarz inequality, we conclude, defining $C_I = C_B^{-1}$,
  \begin{equation*}
    \norm[E]{\Delta w} \leq
    C_I h_{E}^{-1} \sup_{v\in\sobh[0]{1}{E}}
    \frac{\scal[E]{\Delta w}{v}}{\norm[E]{\nabla v}} =
    C_I h_{E}^{-1} \sup_{v\in\sobh[0]{1}{E}}
    \frac{\scal[E]{-\nabla w}{\nabla v}}{\norm[E]{\nabla v}}
    \leq  C_I h_{E}^{-1} \norm[E]{\nabla w} \,.
  \end{equation*}
\end{proof}

\begin{theorem}
  \label{th:interpolationEstim}
  Let $u \in \sobh{s}{\Omega}$, $1\leq s \leq k+1$. Let $\ell_E \in \mathbb{N}$
  be given $\forall E \in \Mh$. There $\exists u_I\in \Vh{k,\ell}$ such that
  \begin{equation}
    \label{eq:interpolationEstim}
    \exists C>0 \colon \norm[\Omega]{u-u_I} + h \norm[\Omega]{\nabla (u-u_I)}
    \leq C h^s \seminorm[s,\Omega]{u} \,.
  \end{equation}
\end{theorem}
\begin{proof}
  Following \cite{Cangiani2017}, let $\Th$ be the sub-triangulation of $\Mh$
  obtained as the union of local sub-triangulations of each polygon $E\in\Mh$,
  linking each vertex to the center of the ball with respect to
  which $E$ is star-shaped. Naturally, $\Th$ inherits the shape-regularity of
  $\Mh$. Let $u_C \in \Poly{k}{\Th}$ be the piecewise polynomial over $\Th$
  defined as the Cl\'ement interpolant of $u$ over $\Th$. This is obtained by
  local projections of $u$ onto piecewise polynomials defined on patches of
  triangles that share a degree of freedom, see \cite{Clement1975}. It holds
  (see \cite[Theorem 1]{Clement1975})
  \begin{equation}
    \label{eq:clementEstim}
    \seminorm[m,\Omega]{u-u_C} \leq C_{Cl,k} h^{s-m} \seminorm[s]{u} \, m\leq s\,,
  \end{equation}
  where $C_{Cl,k}$ depends on the shape-regularity of $\Mh$ and the order
  $k$. Let $w_I\in \sobh{1}{\Omega}$ be the function that solves,
  $\forall E \in \Mh$,
  \begin{equation}
    \label{eq:defwI}
    \begin{cases}
      -\Delta w_I = -\Delta \proj[0,E]{k}{} u_C & \text{in $E$} \,,
      \\
      w_I = u_C & \text{on $\partial E$} \,.
    \end{cases}
  \end{equation}
  By the definition of $w_I$ we have that, $\forall E \in \Mh$,
  $\proj[0,E]{k}{u_C} - w_I$ solves the following problem:
  \begin{equation*}
    \begin{cases}
      -\Delta (\proj[0,E]{k}{}u_C - w_I) = 0 & \text{in $E$} \,,
      \\
      \proj[0,E]{k}{}u_C - w_I = \proj[0,E]{k}{}u_C - u_C & \text{on $\partial E$} \,.
    \end{cases}
  \end{equation*}
  It follows that
  \begin{equation*}
    \norm[E]{\nabla (\proj[0,E]{k}{}u_C - w_I)} =
    \inf\{\norm[E]{\nabla z}\,, z \in \sobh{1}{E}\colon \trace[\partial E]{z}
    = \proj[0,E]{k}{}u_C - u_C \}
    \leq \norm[E]{\nabla (\proj[0,E]{k}{}u_C - u_C)}\,,
  \end{equation*}
  which implies, exploiting the continuity of the operator $\proj[0,E]{k}{}$ and
  \eqref{eq:clementEstim},
  \begin{multline}
    \label{eq:nablauC-wIestim}
    \norm[E]{\nabla(u_C-w_I)} \leq \norm[E]{\nabla(u_C-\proj[0,E]{k}{}u_C)}  +
    \norm[E]{\nabla (\proj[0,E]{k}{}u_C - w_I)}
    \\
    \leq 2 \norm[E]{\nabla\left(\proj[0,E]{k}{}u_C - u_C\right)}
    \leq 2C_{\Pi,k}\norm[E]{\nabla u_C} \leq 2C_{\Pi,k}C_{Cl,k}\norm[E]{\nabla u}\,,
  \end{multline}
  where $C_{\Pi,k}$ depends on the shape-regularity of $E$ and the order $k$.
  We define $u_I \in \Vh{k,\ell}$ such that $\forall E \in \Mh$
  \begin{align}
    \label{eq:uIdef-1}
    \trace[e]{u_I} &= \trace[e]{w_I} \quad \forall e \in \partial E \,,
    \\
    \label{eq:uIdef-2}
    \scal[E]{u_I - w_I}{p} &= 0 \quad \forall p \in \Poly{k-2}{E} \,,
    \\
    \label{eq:uIdef-3}
    \scal[E]{u_I - \proj[\nabla,E]{k}{}w_I}{p} &= 0
                                                 \quad \forall p \in \Poly{k+\ell_E,k-2}{E} \,.
  \end{align}
  Notice that applying Green's theorem we get from \eqref{eq:uIdef-1} and
  \eqref{eq:uIdef-2} that $\proj[\nabla,E]{k}{}u_I = \proj[\nabla,E]{k}{}w_I$ nd
  thus \eqref{eq:uIdef-3} implies
  \begin{equation}
    \label{eq:enhancementuI}
    \scal[E]{u_I - \proj[\nabla,E]{k}{}u_I}{p} = 0 \quad
    \forall p \in \Poly{k+\ell_E,k-2}{E} \,.
  \end{equation}
  We now prove \eqref{eq:interpolationEstim} for $u_I$. Concerning
  $\sobh{1}{\Omega}$-seminorm of $u-u_I$, recalling \eqref{eq:clementEstim} we
  get
  \begin{equation*}
    \begin{split}
      \norm[\Omega]{\nabla (u-u_I)} &\leq \norm[\Omega]{\nabla(u-u_C)} + \norm[\Omega]{\nabla(u_C-u_I)}
      \\
                                    &\leq C_{Cl,k} h^{s-1}\seminorm[s,\Omega]{u} + 
                                      \left(
                                      \sum_{E\in\Mh}\norm[E]{\nabla(u_C-u_I)}^2
                                      \right)^{\frac12}\,.
    \end{split}
  \end{equation*}
  Moreover, from the definition of $w_I$ \eqref{eq:defwI} and the definition of
  $u_I$ \eqref{eq:uIdef-1} we get
  $\trace[e]{u_C} = \trace[e]{w_I} = \trace[e]{u_I}$ for each edge $e$ of
  $\Mh$. Thus, $u_C-u_I\in\sobh[0]{1}{E}$ $\forall E\in\Mh$. We can thus
  estimate the $\lebl{\Omega}$-norm of $u-u_I$ as follows:
  \begin{equation}
    \label{eq:estimu-uI:first}
    \begin{split}
      \norm[\Omega]{u-u_I} &\leq \norm[\Omega]{u-u_C} + \norm[\Omega]{u_C-u_I}
      \\
                           &\leq C_{Cl,k} h^s\seminorm[s,\Omega]{u} + 
                             \left(
                             \sum_{E\in\Mh}\norm[E]{u_C-u_I}^2
                             \right)^{\frac12}
      \\
                           &\leq C_{Cl,k} h^s\seminorm[s,\Omega]{u} + h\left(
                             \sum_{E\in\Mh}C_{p,E}\norm[E]{\nabla(u_C-u_I)}^2
                             \right)^{\frac12}\,,
    \end{split}
  \end{equation}
  where $C_{p,E}$ depends on the shape-regularity of $E$. We now focus on
  estimating $\norm[E]{\nabla(u_C-u_I)}$ $\forall E \in \Mh$. Applying
  \eqref{eq:nablauC-wIestim} we get
  \begin{equation}
    \label{eq:estimuC-uI-1}
    \begin{split}
      \norm[E]{\nabla(u_C-u_I)}
      &\leq \norm[E]{\nabla(u_C-w_I)}
        + \norm[E]{\nabla(w_I-u_I)}
      \\
      & \leq
        2C_{\Pi,k}C_{Cl,k} \norm[E]{\nabla u} + \norm[E]{\nabla(w_I - u_I)} \,.
    \end{split}
  \end{equation}
  Finally, the last term can be bounded applying Green's theorem (recall that
  $w_I-u_I\in\sobh[0]{1}{E}$), \eqref{eq:uIdef-2}, \eqref{eq:enhancementuI}, the
  fact that $\proj[\nabla,E]{k}{}u_I = \proj[\nabla,E]{k}{}w_I$, a
  Cauchy-Schwarz inequality, the approximation properties of polynomial
  projections, \eqref{eq:nablauC-wIestim} and the inverse inequality
  \eqref{eq:inverse-estim} we get
  \begin{equation*}
    \begin{split}
      \norm[E]{\nabla(w_I - u_I)}^2
      &= -\scal[E]{\Delta(w_I-u_I)}{w_I-u_I}
      \\
      &= -\scal[E]{\Delta(w_I-u_I) - \proj[0,E]{k-2}{\Delta(w_I-u_I)}}{w_I-u_I}
      \\
      &= -\scal[E]{\Delta(w_I-u_I) - \proj[0,E]{k-2}{\Delta(w_I-u_I)}}{w_I-\proj[\nabla,E]{k}{}u_I}
      \\
      &= -\scal[E]{\Delta(w_I-u_I) - \proj[0,E]{k-2}{\Delta(w_I-u_I)}}{w_I-\proj[\nabla,E]{k}{}w_I}
      \\
      & \leq \norm[E]{\Delta(w_I-u_I) - \proj[0,E]{k-2}{\Delta(w_I-u_I)}}
        \norm[E]{w_I-\proj[\nabla,E]{k}{}w_I}
      \\
      & \leq C \norm[E]{\Delta(w_I-u_I)} \cdot h_E \norm[E]{\nabla w_I}
      \\
      & \leq C \norm[E]{\Delta(w_I-u_I)} \cdot h_E \norm[E]{\nabla u}
      \\
      & \leq C \norm[E]{\nabla(w_I-u_I)} \norm[E]{\nabla u} \,.
    \end{split}
  \end{equation*}
  Then, collecting \eqref{eq:estimu-uI:first}, \eqref{eq:estimuC-uI-1} and the
  last estimate, we obtain \eqref{eq:interpolationEstim}.
\end{proof}

The interpolation estimate provided by Theorem \ref{th:interpolationEstim} along
with approximation results analogous to the ones obtained in
\cite{BBBPSsupg,Beirao2021} are used to prove the following a priori error
estimates, whose proof is omitted since it is analogous to the one in the cited
references.

\begin{theorem}\label{Th:errorestimates}
  Assume $u\in\sobh{s+1}{\Omega}$ and $f\in\sobh{s-1}{\Omega}$. Then, under the
  current regularity assumptions and if $\ell_E$ is chosen $\forall E\in\Mh$ in
  such a way that \eqref{eq:localcoercivityallk} holds,
  \begin{equation*}
    \ennorm{u-u_h} \leq C h^s \left(\max_{E\in\Th} \{1,\sqrt{h_E\beta_E}\}
    \norm[s+1]{u} + C_{f,\K\beta}\norm[s-1]{f}\right) \,,
  \end{equation*}
  where $C$ is independent of $h$ and on the problem coefficients and
  $C_{f,\K\beta}$ depends on local variations of the problem coefficients.
\end{theorem}


\section{Numerical Results} \label{sec:numres}
In this section we present two numerical experiments, namely Test 1 in Section~\ref{Test1} and Test 2 in Section~\ref{Test2}.
In the first test we confirm the convergence rates predicted by the a priori error analysis of Section~\ref{sec:apriori} and we compare our method with the standard SUPG virtual element discretization \cite{BBBPSsupg} in terms of relative energy errors. In the second test we consider a classic problem taken from \cite{franca1992stabilized} involving approximation of internal and boundary layers. In all the numerical tests our aim is to assess the robustness of the approach in case of advection dominated regime. Therefore, all the benchmark problems are characterized by high mesh P\'eclet numbers. Moreover, for each type of polygon in each mesh considered in the tests, we have done a preliminary assessment of the minimum $\ell_E$ that satisfies Assumption \ref{eq:assumcoercivity}. This is done computing the two smallest eigenvalues of the local stiffness matrix $A^E$, defined as $A^E_{ij}=\scal[E]{\proj[0,E]{k+\ell_E-1}{}\nabla \phi_i}{\proj[0,E]{k+\ell_E-1}{}\nabla \phi_j}$, where $\{\phi_i\}$ denotes the set of local basis functions of the VEM space on a polygon $E\in\Mh$.
We consider \EEVEM  of different orders from one to three and a unit square domain $\Omega =  (0,1) \times (0,1)$.

\subsection{Test 1}\label{Test1}
\begin{figure}[H]
\centering
\includegraphics[width=10cm]{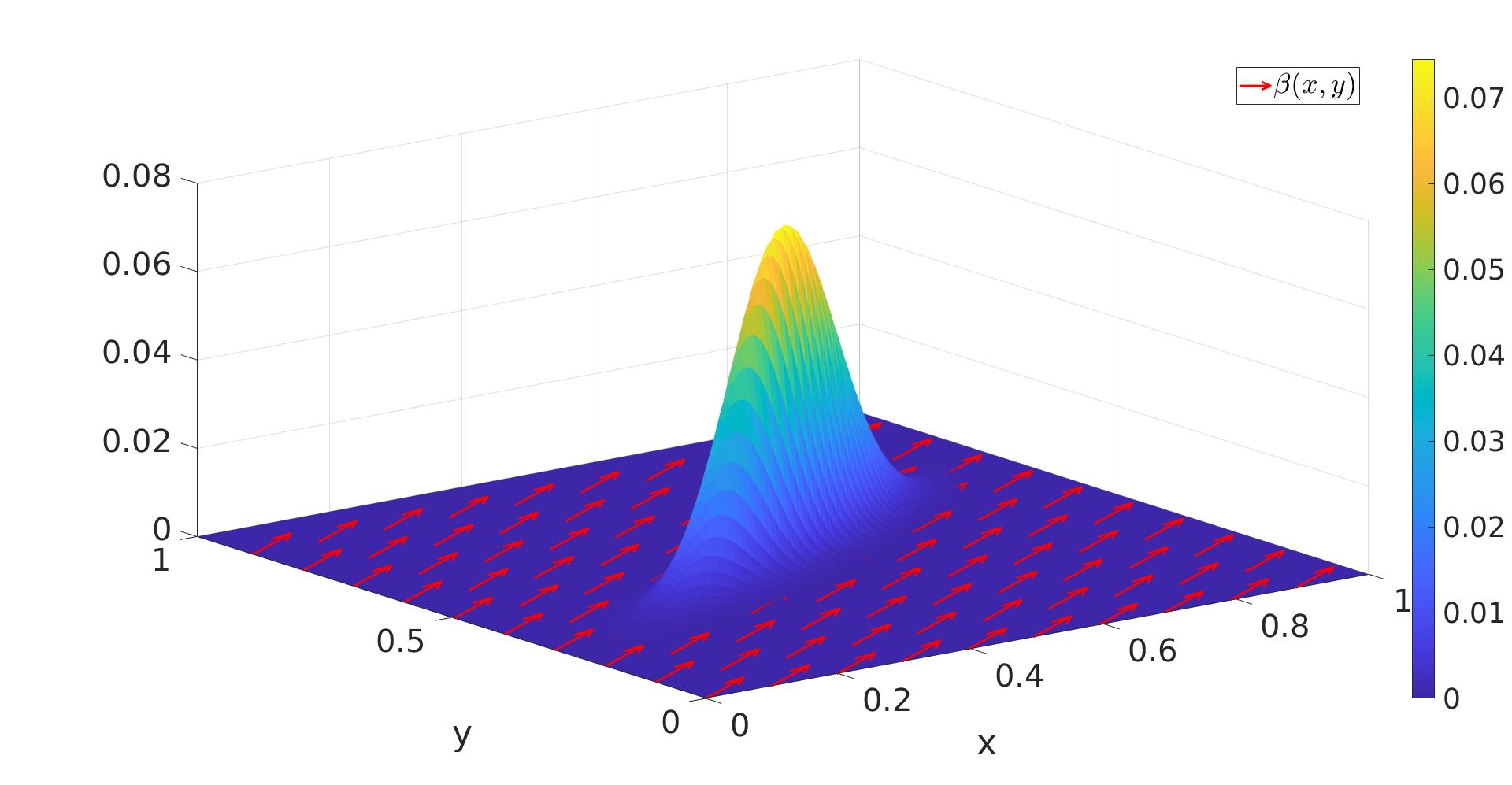}
\caption{Test 1: exact solution $u(x,y)$ and transport velocity field $\beta(x,y)$ (red arrows).}\label{Fig1}
\end{figure}

In this first test, we consider an advection-diffusion problem characterized by a diffusivity parameter $\K = 10^{-9}$ and a transport velocity field $\beta = (1, 0.545)$. The size of the meshes is chosen to guarantee that for the selected value of $\K$ the mesh P\'eclet number is much greater than one for all $k$. 
The forcing term $f$ and the boundary conditions are such that the exact solution (depicted in Figure~\ref{Fig1}) is
\begin{equation*}
u(x,y) = c_1 \ x y \ (x-1)(y-1) \ e^{-c_2 (c_4(c_2-x)^2+c_3(c_2-y)^2-c_3(c_2-x)(c_2-y))},
\end{equation*}
where $c_1 = \frac{3}{\sqrt{2 \pi}}$, $c_2 = \frac{1}{2}$, $c_3 = 1000$ and $c_4 = \frac{1}{3.3} \cdot 10^3$.
The selected solution $u$ is characterized by a strong anisotropy. Indeed, in Figure~\ref{Fig1}, we can see that the solution exhibits a strong boundary layer in a direction approximately perpendicular to the direction of the transport velocity field $\beta$.

\begin{figure}[H]
\centering
\captionsetup{justification=centering}
\subfloat[][$\mathcal{T}_1$]{
\includegraphics[width=4cm]{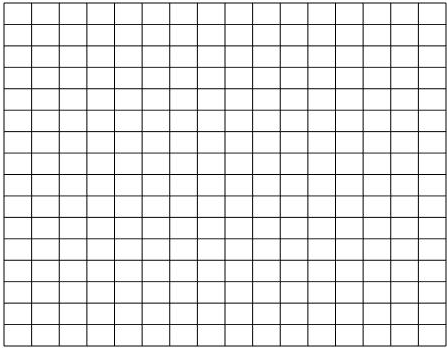}\label{Fig2_1}}
\subfloat[][$\mathcal{T}_2$]{
\includegraphics[width=4cm]{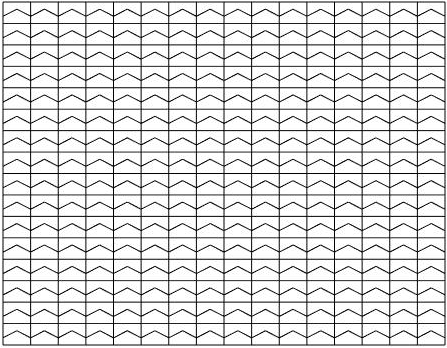}\label{Fig2_2}}
\subfloat[][$\mathcal{T}_3$]{
\includegraphics[width=4cm]{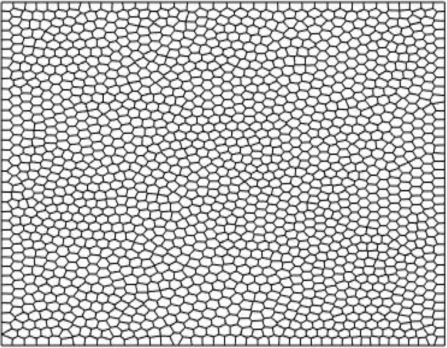}\label{Fig2_3}}
\caption{Meshes.}
\label{Fig2}
\end{figure}

To study the convergence of the method, we consider three different families of meshes ($\mathcal{T}_{1}$, $\mathcal{T}_2$ and $\mathcal{T}_3$) and four different refinements for each one of them. The first mesh of each sequence is reported in Figure~\ref{Fig2}.

\begin{table}[H]
\footnotesize
\begin{center}
\begin{threeparttable}[b]
\centering
\begin{tabular}{ccccccccccc}
\toprule
\multicolumn{1}{c}{} 
& \multicolumn{3}{c}{\textbf{$\mathcal{T}_1$}} 
& \multicolumn{3}{c}{\textbf{$\mathcal{T}_2$}} 
& \multicolumn{3}{c}{\textbf{$\mathcal{T}_3$}} \\
  \cmidrule(rl){2-4} \cmidrule(rl){5-7} \cmidrule(rl){8-10}
\textbf{$k$}  &     1 &  2 &  3-4 &   1 &  2 &  3-4 &   1 &  2 &  3-4\\
\midrule
$\mathcal{M}_{first}$  & $2 \cdot 10^{7}$ & $4 \cdot 10^{6}$ & $1 \cdot 10^{6}$ & $1 \cdot 10^{7}$ & $3 \cdot 10^{6}$ & $8 \cdot 10^{5}$ & $6 \cdot 10^{6}$ & $2 \cdot 10^{6}$ & $4 \cdot 10^{5}$ \\
$\mathcal{M}_{last}$  &   $2 \cdot 10^{6}$ & $5 \cdot 10^{5}$&  $1 \cdot 10^{5}$ &   $2 \cdot 10^{6}$ & $4 \cdot 10^{5}$&  $1 \cdot 10^{5}$  & $2 \cdot 10^{6}$ &  $5 \cdot 10^{5}$ & $1 \cdot 10^{5}$\\
\bottomrule
\end{tabular}
\begin{tablenotes}
     \end{tablenotes}
\end{threeparttable}
\caption{Test 1: Mean values of the mesh P\'eclet number for the first mesh $\mathcal{M}_{first}$ and last mesh $\mathcal{M}_{last}$ of the mesh families $\mathcal{T}_1$, $\mathcal{T}_2$ and $\mathcal{T}_3$.}\label{P_table}
\end{center}
\end{table}

The mesh family $\mathcal{T}_1$ consists of standard cartesian elements, the mesh family $\mathcal{T}_2$ is composed of both concave and convex polygons and the mesh family $\mathcal{T}_3$ has been constructed using Polymesher \cite{Polymesher}. The first two groups of meshes are refined splitting the existing elements in half. For the last group of meshes this approach is not feasible. Therefore, we simply increase the number of elements by means of Polymesher. Consequently, the tessellations of this family include polygons having different numbers of edges.
 Moreover, in Table~\ref{P_table} we report also the mean mesh P\'eclet number for the first mesh and the last mesh of each mesh family.

\begin{table}[H]
\begin{center}
\begin{threeparttable}[b]
\centering
\begin{tabular}{ccccccccccc}
\toprule
\multicolumn{2}{c}{} 
& \multicolumn{2}{c}{\textbf{$\mathcal{T}_1$}} 
& \multicolumn{2}{c}{\textbf{$\mathcal{T}_2$}} 
& \multicolumn{4}{c}{\textbf{$\mathcal{T}_3$}} \\ 
  \cmidrule(rl){7-10}
\textbf{$N^{V}_{E}$}  &  & 4 &  & 5  &  &1-4 &5 & 6 &7\\
\midrule
$k = 1$ & & 1 &  & 1 &  & 1 & 1 & 2 & 2 \\
$k = 2$ & & 2 &  & 1 &  & 1 & 1 & 2 & 2 \\
$k = 3$ & & 2 &  & 1 &  & 1 & 1 & 2 & 2 \\
$k = 4$ & & 2 &  & 2 &  & 1 & 2 & 3 & 4 \\
\bottomrule
\end{tabular}
\end{threeparttable}
\caption{Values of $\ell_E$ for the proposed method and tessellations $\mathcal{T}_1$, $\mathcal{T}_2$ and $\mathcal{T}_3$, related to the number of vertices $N_{E}^{V}$ of the polygons in each mesh.}\label{l_table}
\end{center}
\end{table}

Now, for each type of polygon of the meshes analyzed, we have to discuss the proper choice of $\ell_E$.
This choice has been done, as said previously, selecting the value $\ell_E$ that guarantees the numerical coercivity of the local stiffness matrix. 
In Table \ref{l_table}, we report the values of $\ell_E$ computed following this criterion and adopted specifically to solve the problem described in the present subsection. Notice that $\ell_E$ does not depend only on the number of vertices of the polygon, but also on its geometry. Indeed, if we consider the line corresponding to $k=2$ in Table \ref{l_table}, we can see that we require $\ell_E=2$ for $\mathcal{T}_1$, where quadrilaterals are all squares, and $\ell_E=1$ for $\mathcal{T}_3$, that features generally shaped quadrilaterals.

\begin{figure}[H]
\centering
\captionsetup{justification=centering}
\subfloat[][$k=1$]{
\includegraphics[width=0.45\linewidth]{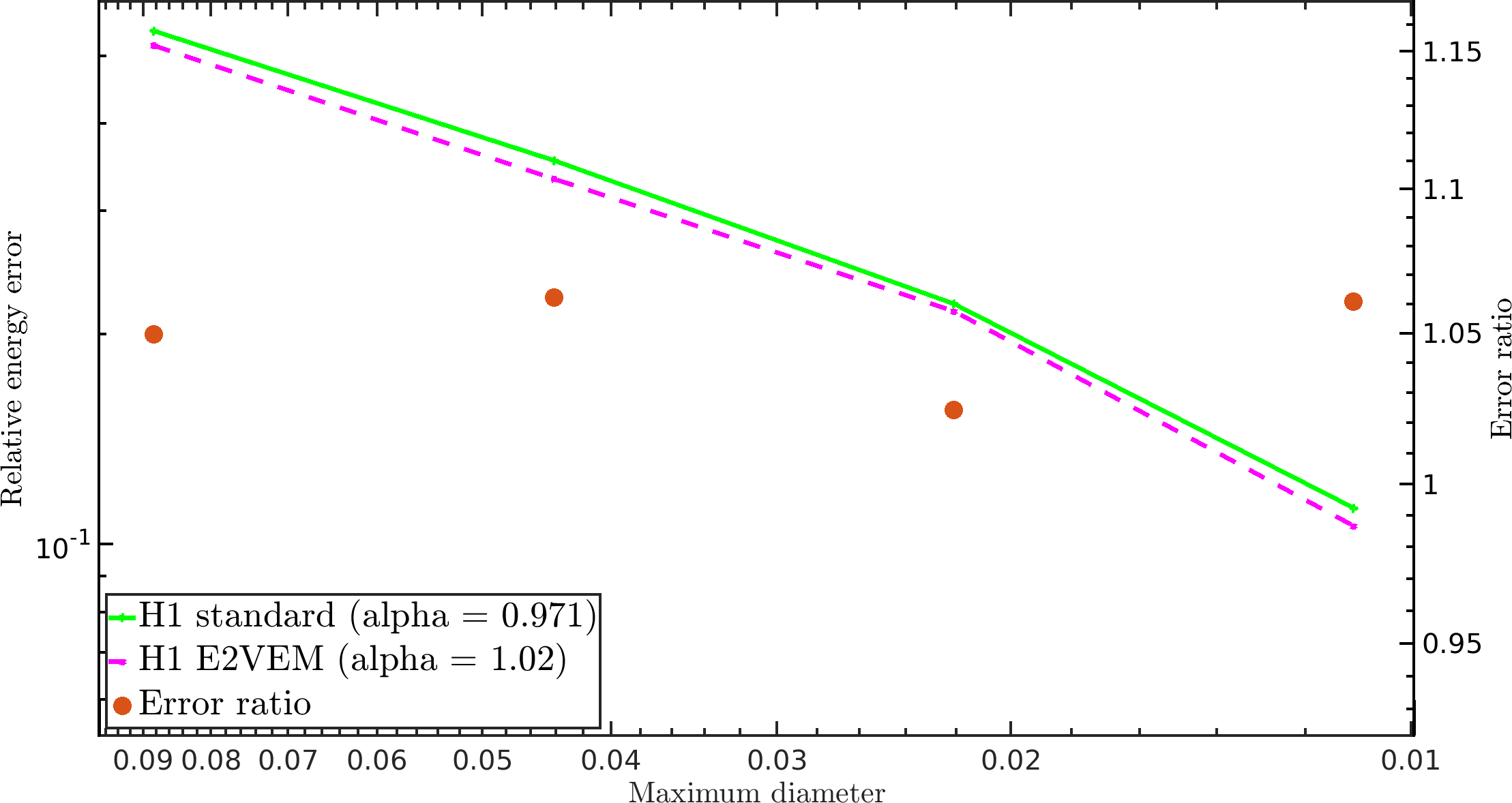}\label{Fig3_1}}\hfill
\subfloat[][$k=2$]{
\includegraphics[width=0.45\linewidth]{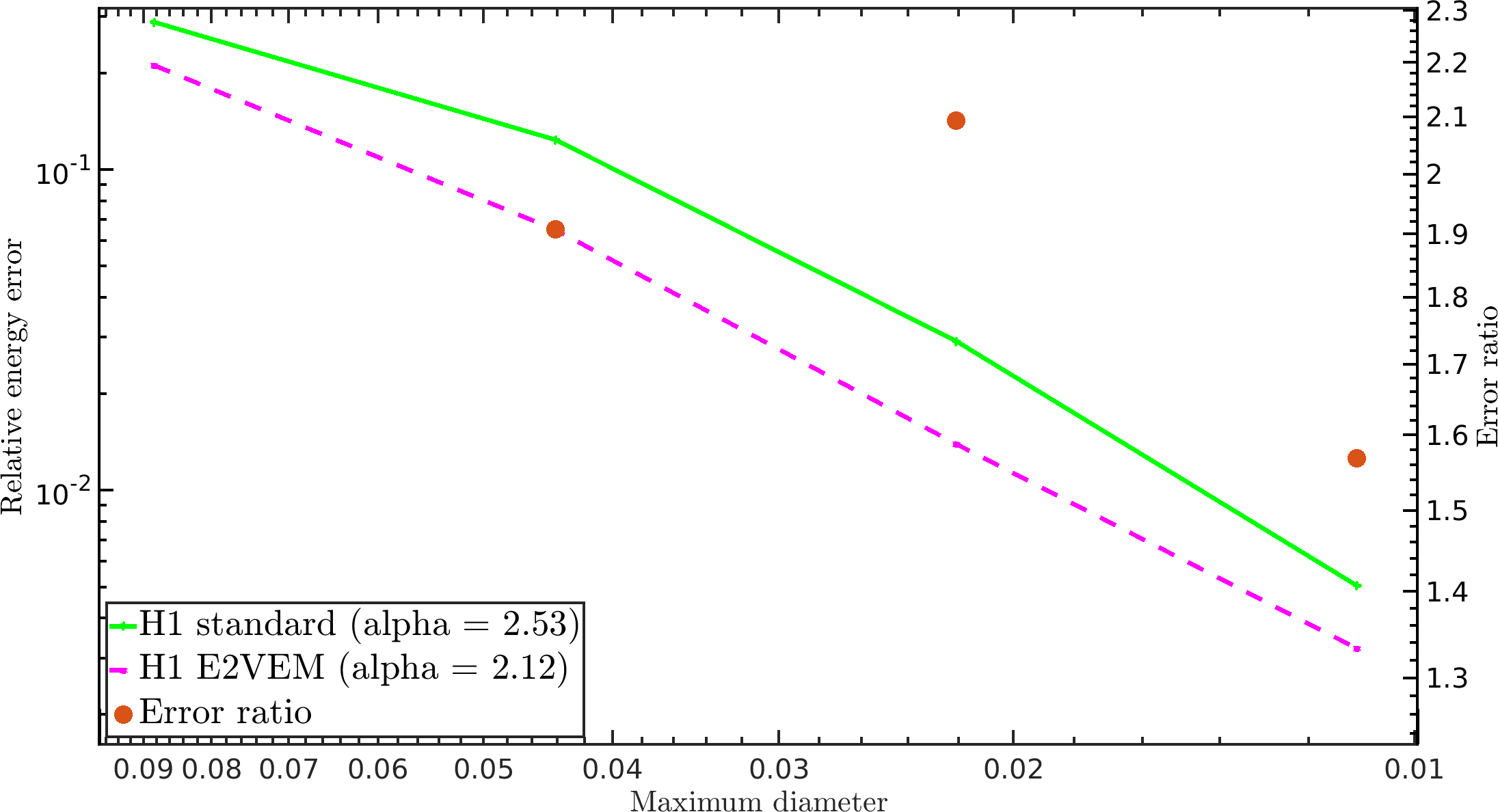}\label{Fig3_2}}
\qquad
\subfloat[][$k=3$]{
\includegraphics[width=0.45\linewidth]{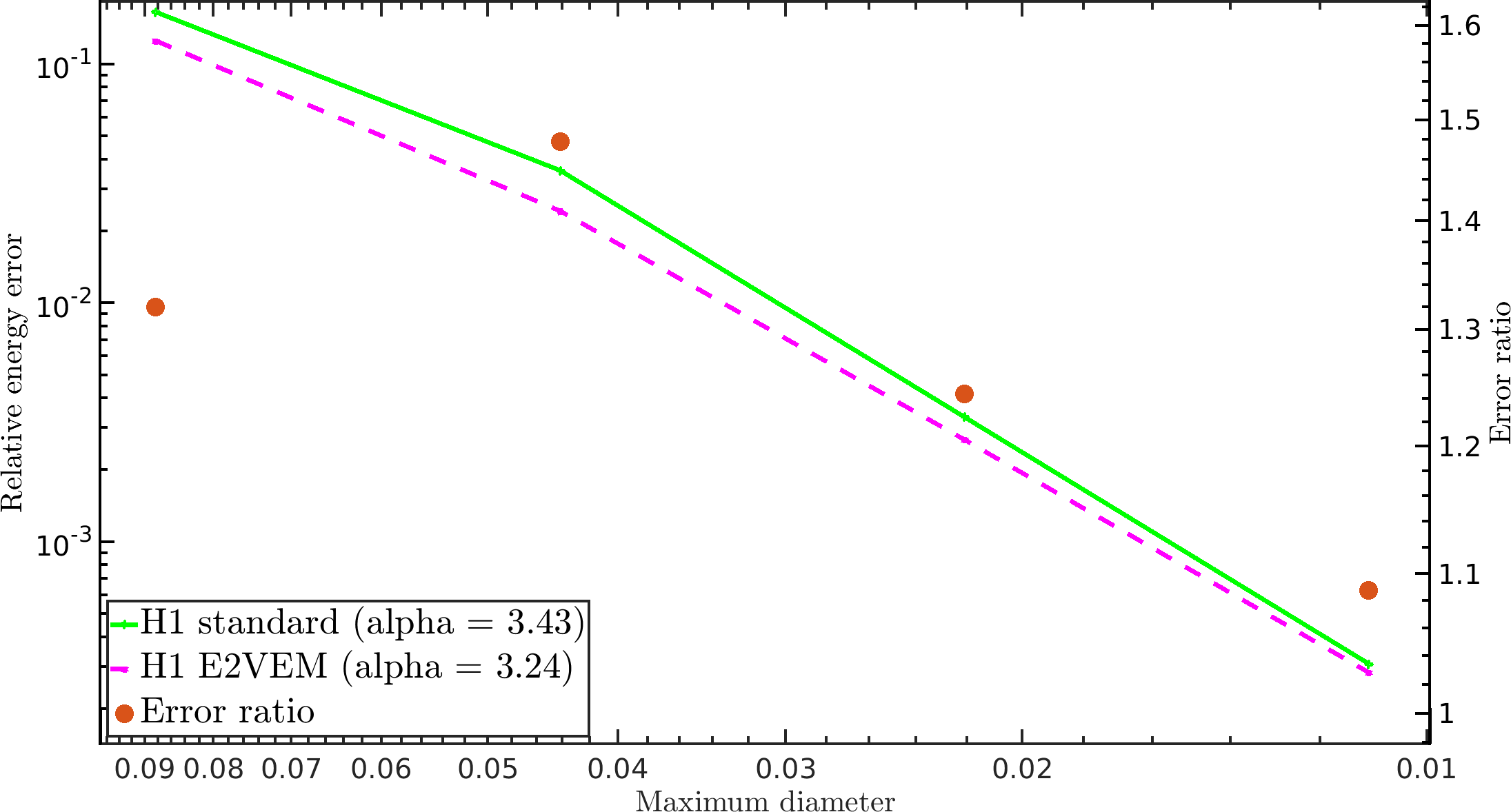}\label{Fig3_3}}\hfill
\subfloat[][$k=4$]{
\includegraphics[width=0.45\linewidth]{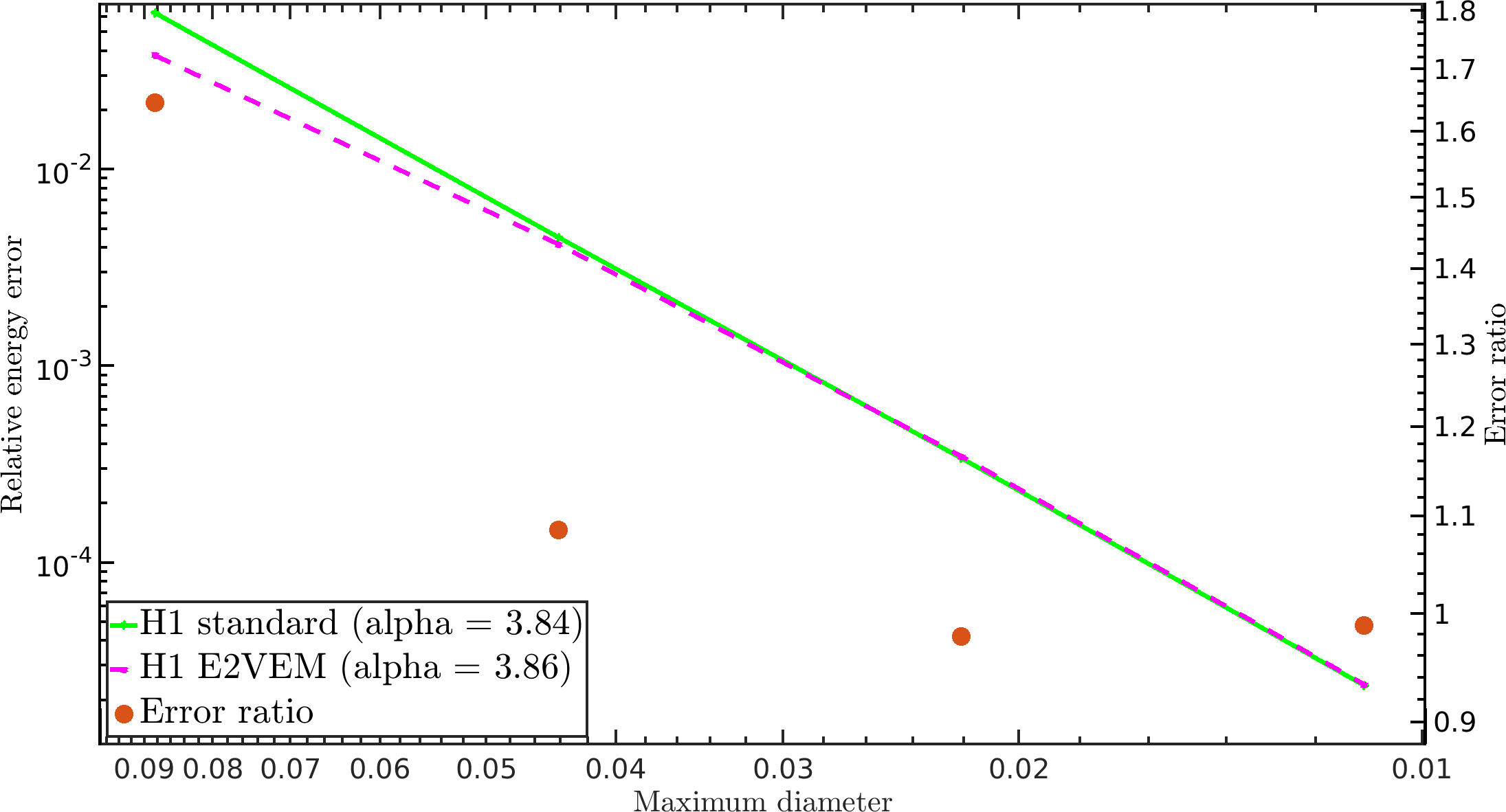}\label{Fig3_4}}
\caption{Test 1: convergence curves (tessellation $\mathcal{T}_1$).}
\label{Fig3}
\end{figure}

\begin{figure}[H]
\centering
\captionsetup{justification=centering}
\subfloat[][$k=1$]{
\includegraphics[width=0.45\linewidth]{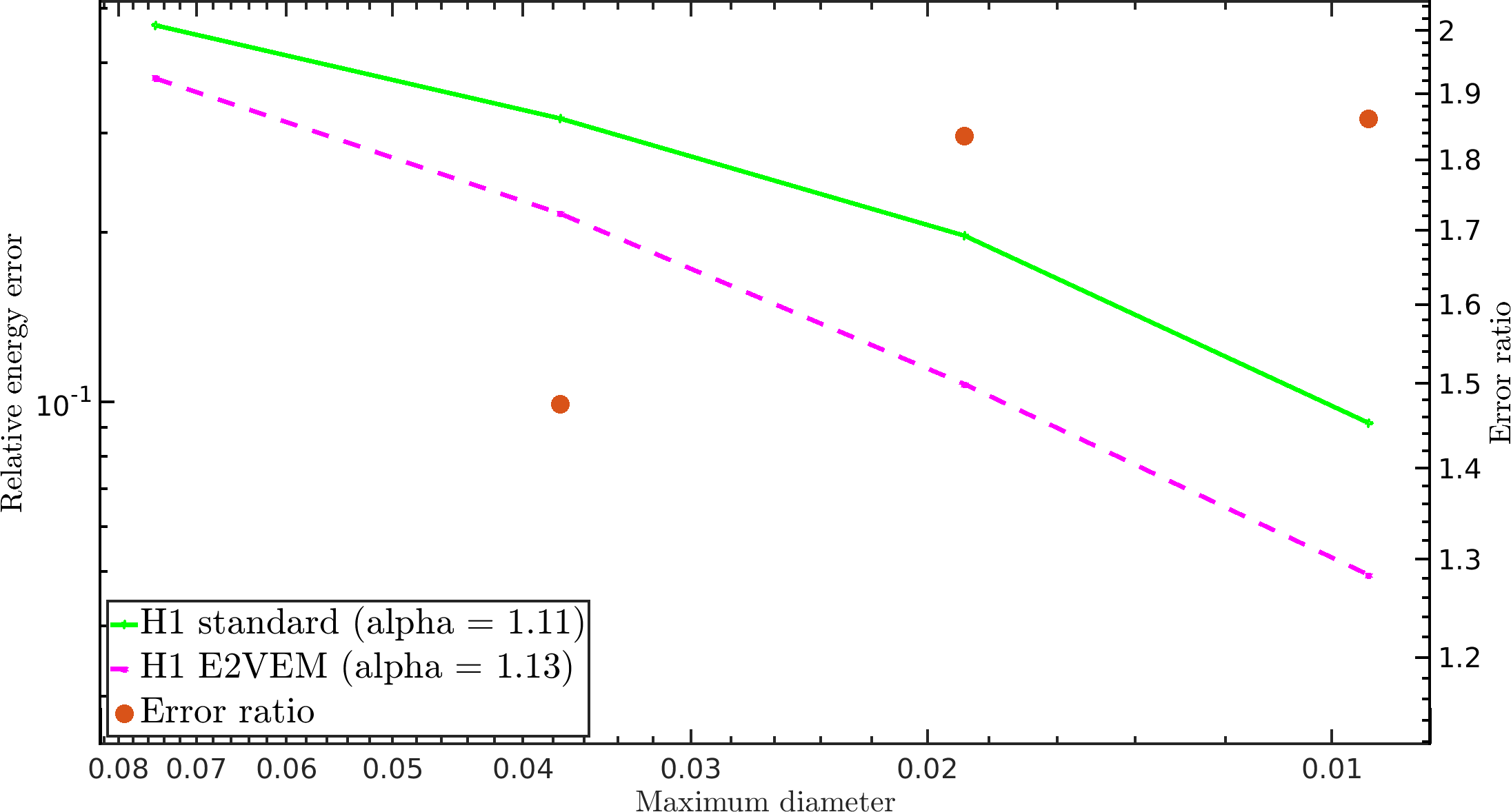}\label{Fig5_1}}\hfill
\subfloat[][$k=2$]{
\includegraphics[width=0.45\linewidth]{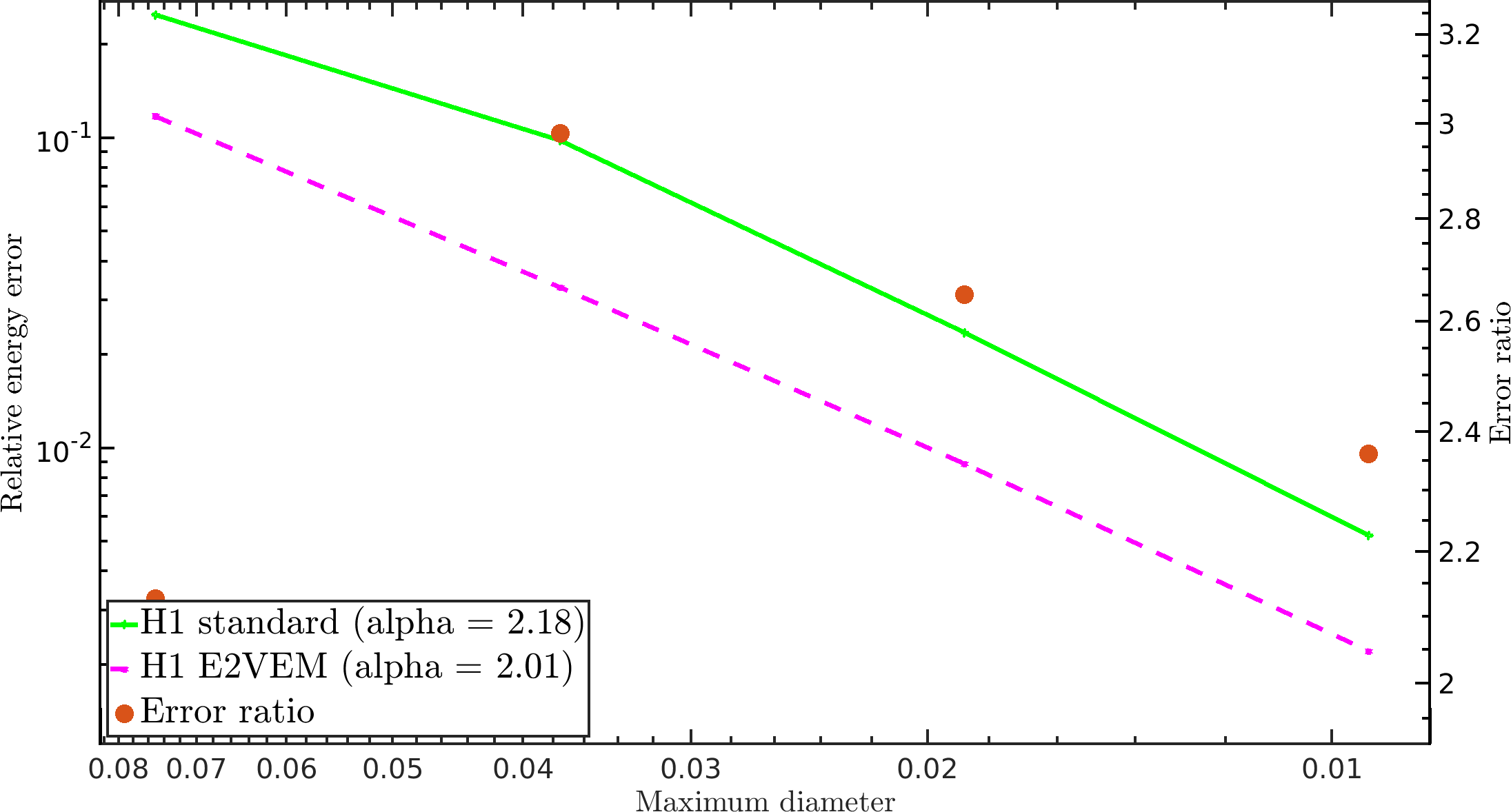}\label{Fig5_2}}
\qquad
\subfloat[][$k=3$]{
\includegraphics[width=0.45\linewidth]{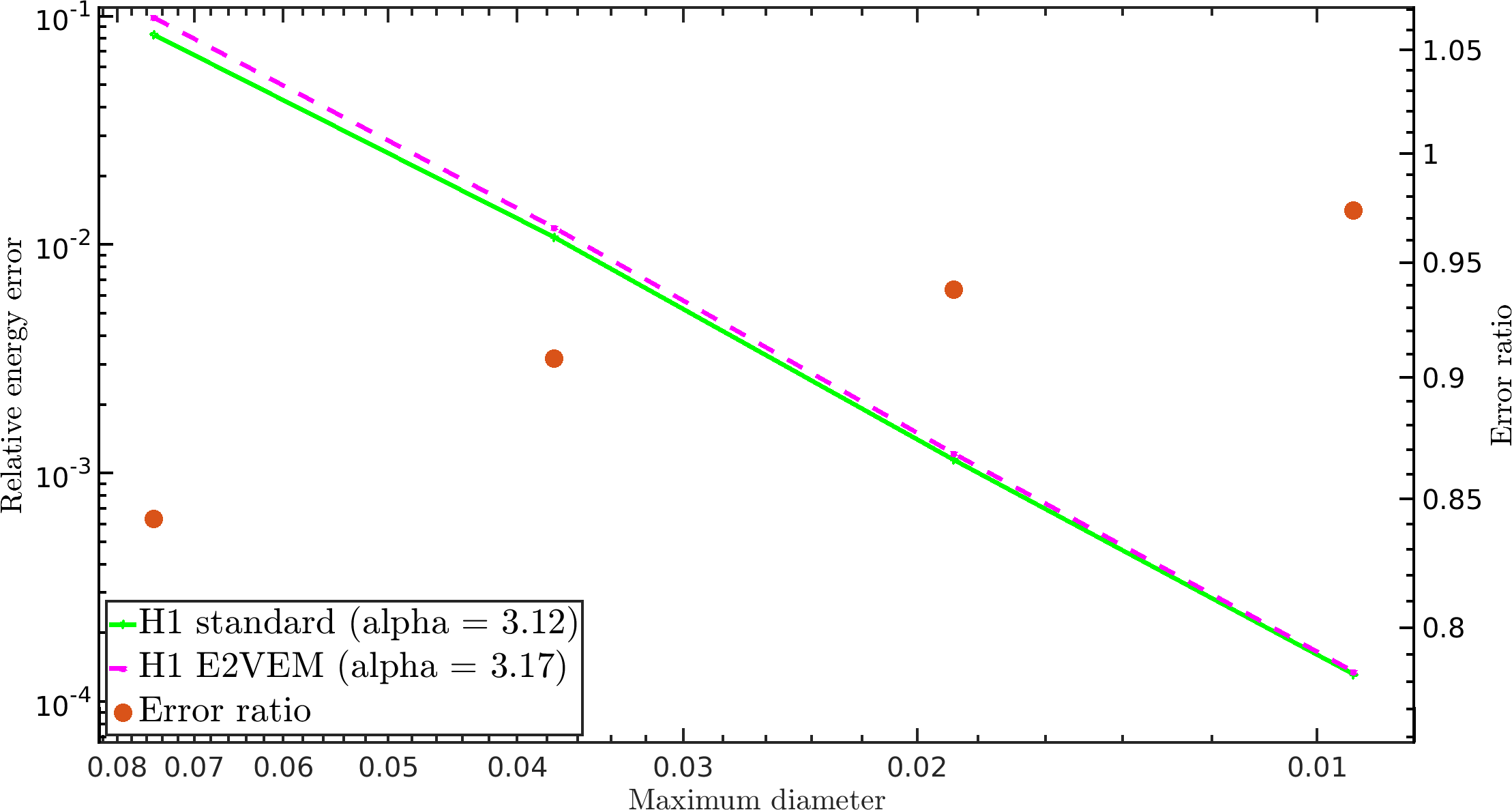}\label{Fig5_3}}\hfill
\subfloat[][$k=4$]{
\includegraphics[width=0.45\linewidth]{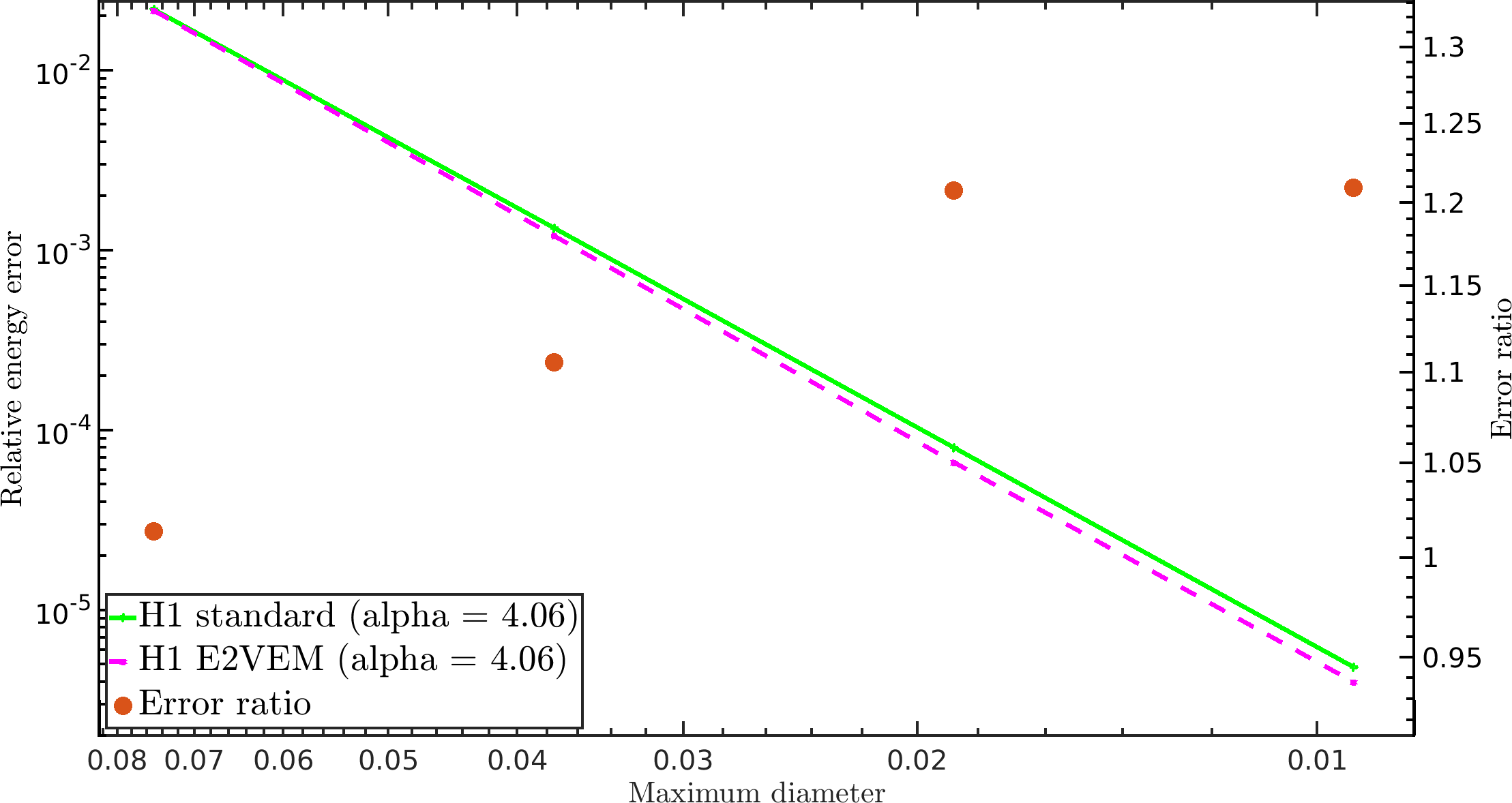}\label{Fig5_4}}
\caption{Test 1: convergence curves (tessellation $\mathcal{T}_2$).}
\label{Fig5}
\end{figure}

\begin{figure}[H]
\centering
\captionsetup{justification=centering}
\subfloat[][$k=1$]{
\includegraphics[width=0.45\linewidth]{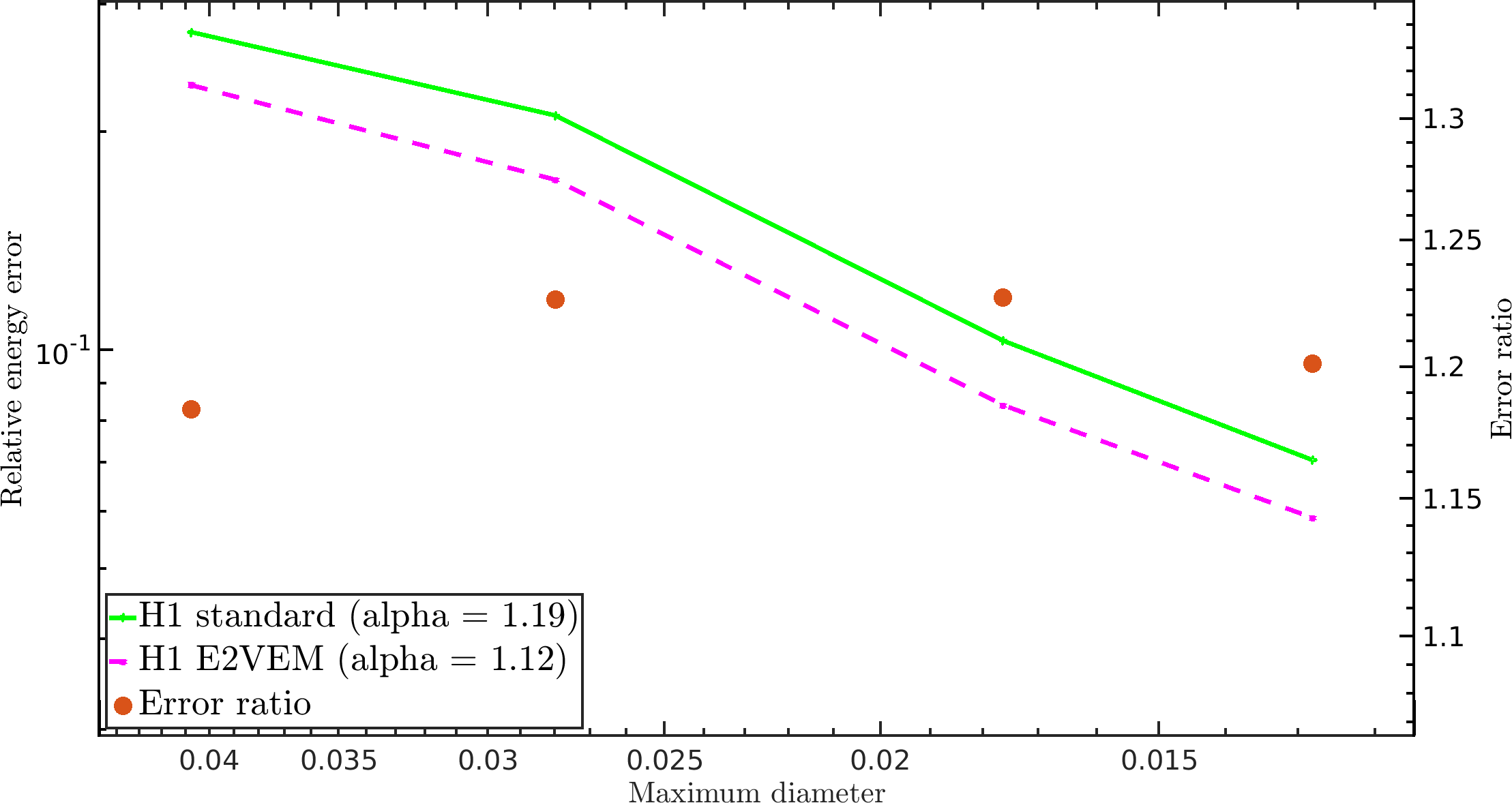}\label{Fig4_1}}\hfill
\subfloat[][$k=2$]{
\includegraphics[width=0.45\linewidth]{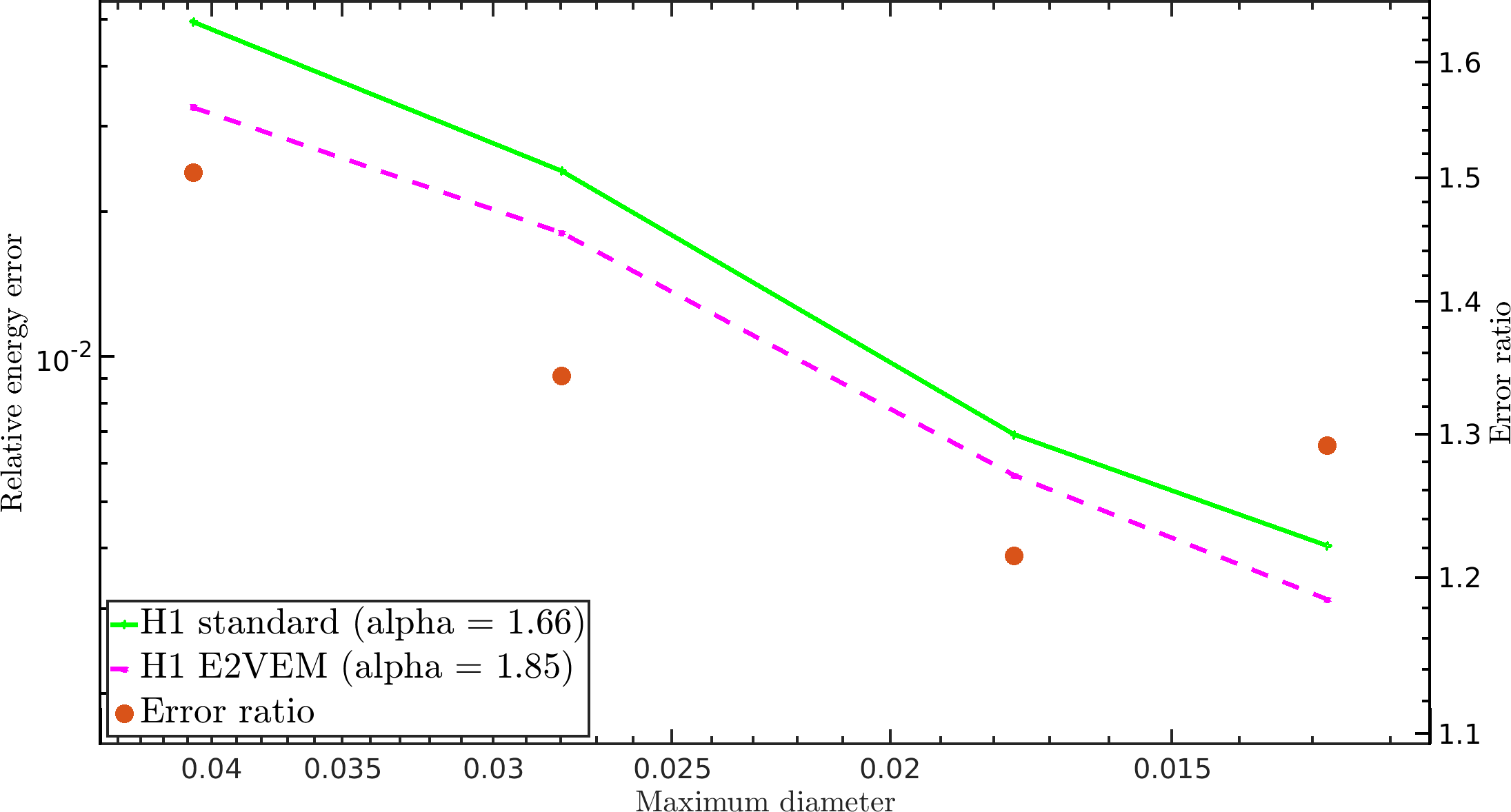}\label{Fig4_2}}
\qquad
\subfloat[][$k=3$]{
\includegraphics[width=0.45\linewidth]{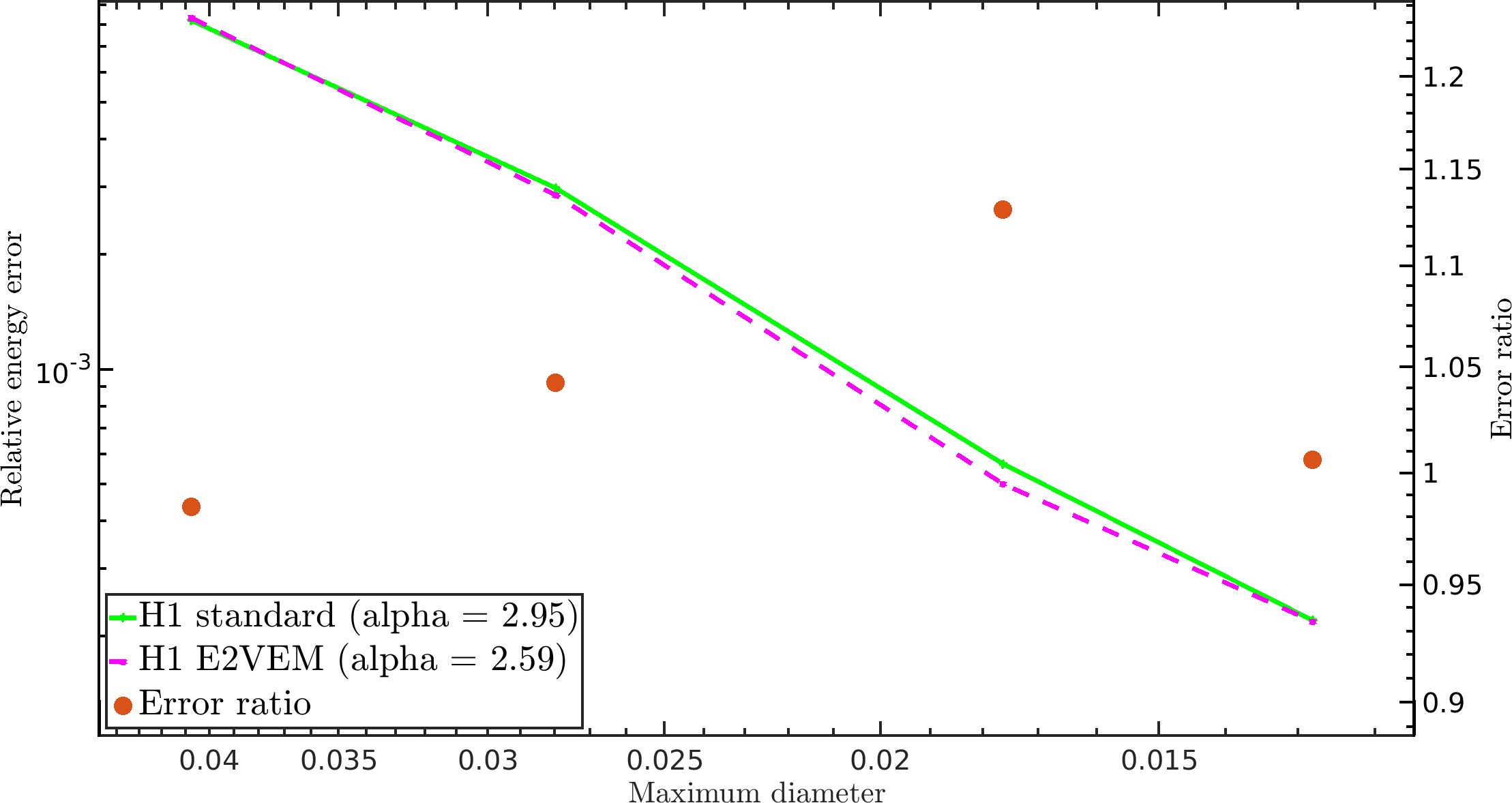}\label{Fig4_3}}\hfill
\subfloat[][$k=4$]{
\includegraphics[width=0.45\linewidth]{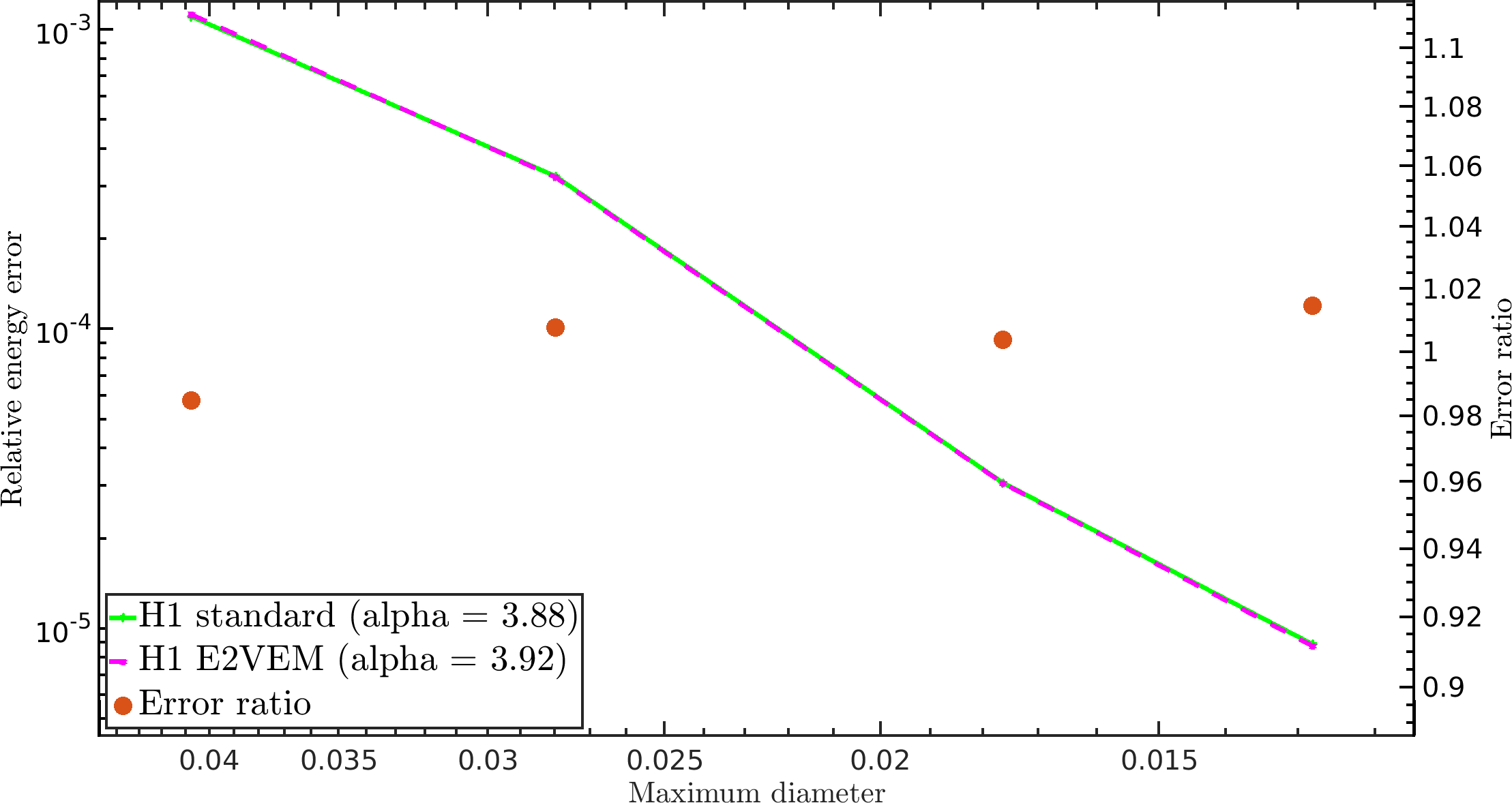}\label{Fig4_4}}
\caption{Test 1: convergence curves (tessellation $\mathcal{T}_3$).}
\label{Fig4}
\end{figure}

In Figures~\ref{Fig3}, \ref{Fig5} and \ref{Fig4} we show the convergence curves in log-log scale. We report the relative errors computed in the energy norm plotted against the maximum diameter of the discretization for both the \EEVEM and the classical VEM \cite{BBBPSsupg} and for orders from $1$ to $4$. The computed relative error is based on the difference between the exact solution and the projection $\proj[\nabla,E]{k}{}$ of the discrete solution $u_{h}$ and it is given by  the following expression
\begin{equation*}
  err = \displaystyle
  \sqrt{\frac{ \sum_{E\in\Mh}
      \norm[E]{\sqrt{\K} \nabla (u - \proj[\nabla,E]{k}{} u_{h})}^2 +
      \tau_{E} \norm[E]{ \beta \cdot \nabla (u - \proj[\nabla,E]{k}{} u_{h})}^{2}}
    { \sum_{E \in \Mh} \norm[E]{\sqrt{\K} \nabla u}^2 +
      \tau_{E} \norm[E]{\beta \cdot \nabla u}^{2}}} \,.
\end{equation*}

The charts of Figures ~\ref{Fig3}, \ref{Fig5} and \ref{Fig4} display two y-axis: the one on the left is related to the relative energy error, whereas the one on the right is related to the ratio between the VEM error and the \EEVEM error. In the legend, \textit{alpha} denotes the numerical rates of convergence, computed using the last two meshes of each refinement. The numerical rates of convergence for both the two methods are in agreement with the theoretical findings for the energy norm of the problem (Theorem \ref{Th:errorestimates}).

Figures~\ref{Fig3_1}, \ref{Fig4_1} and \ref{Fig5_1} show that for order $k=1$ the two methods
provide very close results on the mesh family $\mathcal{T}_1$, whereas the \EEVEM over performs the
classical VEM on the mesh families $\mathcal{T}_2$ and $\mathcal{T}_3$. Thus, the results suggest
that the \EEVEM is able to decrease the magnitude of the error with respect to the classical VEM
when dealing with solutions characterized by strong anisotropies. Indeed, the possibility to avoid
an arbitrary stabilizing part enhances the performance of the method. An analogous trend is observed
for $k=2$, though in this case we observe better performances also for the cartesian mesh
$\mathcal{T}_1$. Moreover, in Figures~\ref{Fig5_1} and \ref{Fig5_2} we notice that the error
difference between the \EEVEM and the VEM seems to be even stronger than the one observed in
Figures~\ref{Fig3_1}, \ref{Fig3_2}, \ref{Fig4_1} and \ref{Fig4_2}. We believe that this could be
related to the presence of non-convex elements.

 Figures~\ref{Fig3_3}, \ref{Fig3_4}, \ref{Fig5_3}, \ref{Fig5_4}, \ref{Fig4_3} and \ref{Fig4_4} show that for $k = 3$ and $k = 4$ the VEM and the \EEVEM exhibit an almost equivalent behaviour for all the considered meshes. A similar trend was also observed for anisotropic elliptic problems in \cite{BBME2VEMLetter}. Indeed, for higher orders we expect that the polynomial part of the VEM bilinear form tends to predominate the stabilizing part reducing the error gap between the two methods.

\subsection{Test 2}\label{Test2}
For the second test, we consider a classic benchmark problem in the SUPG literature characterized by the presence of different layers. This problem was originally proposed in the context of standard Finite Element Methods in \cite{franca1992stabilized}.

\begin{figure}[H]
\centering
\includegraphics[width=7cm]{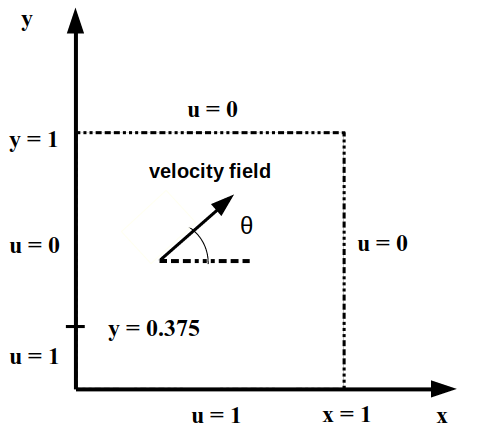}
\caption{Test 2: computational domain and boundary conditions.}
\label{Fig6}
\end{figure}

The computational domain $\Omega$ as well as the boundary conditions are illustrated in Figure~\ref{Fig6}. Notice the discontinuous Dirichlet boundary condition on the left side of the domain. The diffusion coefficient is set equal to $\K = 10^{-6}$, the velocity field is $\beta = (\cos{\theta}, \sin{\theta})$, where $\theta = \frac{\pi}{4}$, and the forcing term $f$ is null. The resulting P\'eclet number is very large (around $10^6$). At the inflow boundary, the velocity propagates the non-homogeneous boundary condition originating an internal boundary layer, whereas at the outflow boundary an outflow boundary layer is produced due to the homogeneous boundary conditions.

\begin{figure}[H]
\centering
\captionsetup{justification=centering}
\subfloat[][$k=1$]{
\includegraphics[width=6.5cm]{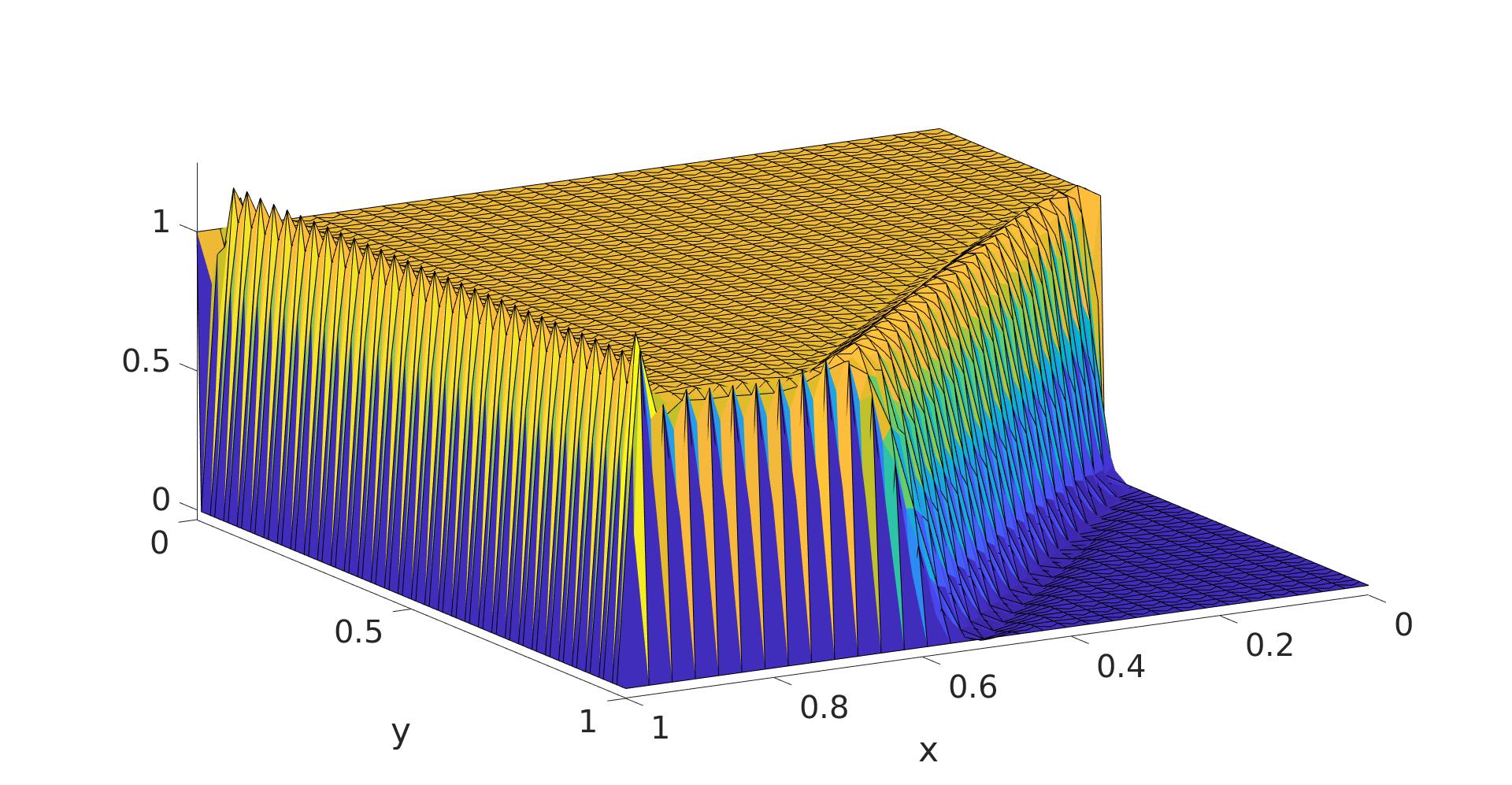}\label{Fig7_1}}
\subfloat[][$k=3$]{
\includegraphics[width=6.5cm]{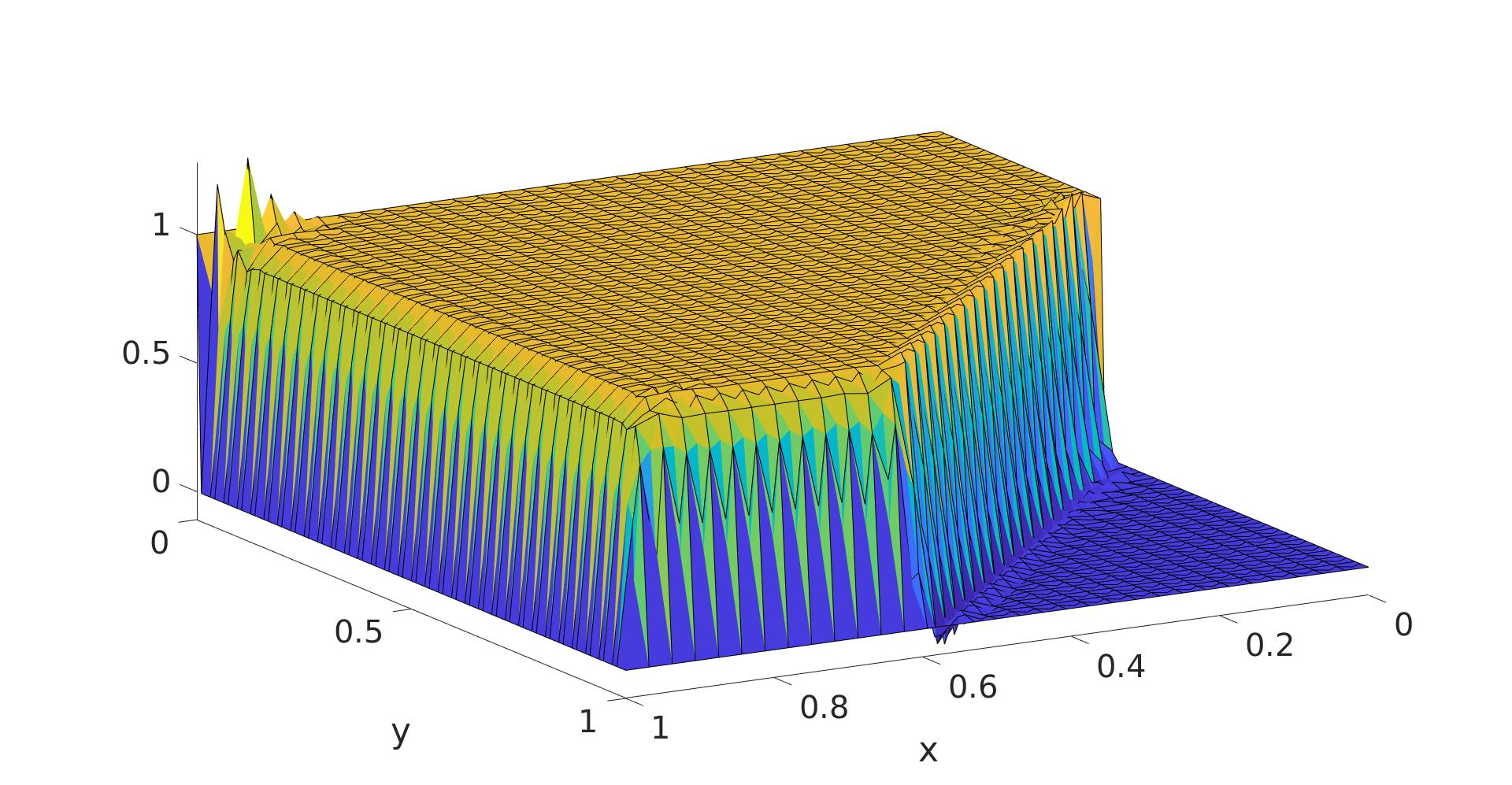}\label{Fig7_3}}
\caption{Test 2: \EEVEM  solution for $k=1$ and $k = 3$ on $\mathcal{T}_2$.}
\label{Fig7}
\end{figure}

We solve the problem using $\mathcal{T}_2$ tessellation, depicted in Figure~\ref{Fig2_2}. Figures~\ref{Fig7_1} and \ref{Fig7_3} show the results obtained for \EEVEM of order $k=1$ and $k=3$. The numerical solutions are comparable to the ones presented in \cite{franca1992stabilized,BBBPSsupg,BBM}. As it is typical for this problem we notice the presence of undershoots and overshoots near the internal boundary layer. Moreover, smoother solutions are obtained increasing the order of the method.

\section{Conclusions}

We presented a numerical investigation of the performances of SUPG-stabilized
stabilization-free VEM, assessing the possibility of using higher-order
polynomial projectors in the definition of the discrete bilinear form and avoiding
the use of a stabilizing bilinear form. We also provided an interpolation
estimate for the new scheme, analogous to the one obtained for standard
VEM. Numerical results show that the possibility of avoiding a stabilizing bilinear form,
that adds some artificial diffusion to the problem in order to ensure coercivity,
can enhance the performances of classical VEM methods when there is the need to
approximate strong boundary layers, as in the case of convection dominated
problems.


\section*{Acknowledgments}
The authors are members of the INdAM-GNCS. A.B and F.M. kindly acknowledge
financial support provided by PNRR M4C2 project of CN00000013 National Centre
for HPC, Big Data and Quantum Computing (HPC) CUP:E13C22000990001.  A.B. kindly
acknowledges partial financial support provided by INdAM-GNCS Projects 2022 and
2023, MIUR project ``Dipartimenti di Eccellenza'' Programme (2018–2022)
CUP:E11G18000350001 and by the PRIN 2020 project (No. 20204LN5N5\_003).


\printbibliography
%

\end{document}